\documentclass[12pt]{article}
\usepackage[utf8]{inputenc}
\usepackage[T1]{fontenc}
\usepackage[leqno]{mathtools}
\usepackage{epic}
\usepackage{fourier}
\usepackage{mathrsfs}
\usepackage{amssymb}
\usepackage{mathrsfs}
\usepackage{graphicx}
\usepackage{a4wide} 
\usepackage{latexsym}
\usepackage[overload]{empheq}
\usepackage{pst-all}
\usepackage[english,francais]{babel}
\newenvironment{preuve}[1]{\par\noindent\textsl
{Preuve #1} \par\noindent}%
{\unskip
\nobreak\hfil\penalty50\hskip2em\null\nobreak\hfil%
$\blacksquare$\parfillskip0pt\par\medskip}

\newtheorem{theorem}{Théorème}
\newtheorem{lem}{Lemme}
\newtheorem{coro}{Corollaire}
\newtheorem{definition}{Définition}
\newtheorem{proposition}{Proposition}
\newcommand{\tore}{\mathbb T}
\newcommand{\cc}{\mathbb C}
\newcommand{\nn}{\mathbb N}
\newcommand{\rr}{\mathbb R}
\newcommand{\zz}{\mathbb Z}
\newcommand{\trace}{\mathrm{tr}}
\newcommand{\myfrac}[2]{\frac{#1}{\rule{0pt}{7pt}#2}}
\newcommand{\spec}{\mathrm{Spec}}

\newcommand{\C}{\mathcal {C}}
\newcommand{\Po}{\mathcal {P}}
\newcommand{\ovl}{\overline}
\newcommand{\So}{\mathbb {S}}
\newcommand{\Lam}{\Lambda}

\newcommand{\la}{\langle}
\newcommand{\ra}{\rangle}
\newcommand{\re}{\mathrm{Re}}\newcommand{\ma}{\mathcal}

\renewcommand{\le}{\leqslant}
\renewcommand{\ge}{\geqslant}
\DeclareMathOperator{\sgn}{\mathrm{signe}}
\DeclareMathOperator{\im}{Im}
\author{\hbox{Jean-Marc  Rinkel} \\{Universit\'{e} de Paris Sud}\\
 { e-mail:jean-marc.rinkel@math.u-psud.fr}
  \and {Abdellatif Seghier }\\{Universit\'{e} de Paris Sud}\\
{  e-mail:Abdellatif.Seghier@math.u-psud.fr}}
 \title{Matrices de Toeplitz tronqu\'ees sur des polygones convexes. Cas du triangle.}
 \date{}
 
\begin{document}
 \maketitle
 \begin{abstract}
\noindent On considère la classe des fonctions positives, essentiellement bornées et semi-continues inférieurement définies sur le tore $\tore^2$. Si $f$ est dans cette classe, elle appara\^{i}t comme le symbole d'un opérateur de Toeplitz tronqué sur un triangle $\Lam$ paramétré par un entier $\lambda$, noté $T_{\Lam}(f)$. Après avoir établi une structure géométrique de l'inverse de cet opérateur, nous donnons un développement asymptotique de la trace de $T_{\Lam}(f)^{-1}$ qui met en évidence la géométrie du triangle. La base de ce théorème est la possibilité de factoriser le symbole $f$ sous la forme $f=\alpha\overline\alpha$ où le spectre de $\alpha$ peut \^{e}tre localisé dans un demi cône donné a priori, autre résultat établi dans cet article. Ce théorème de trace permet en particulier de retrouver le théorème de Linnik-Szegö sur l'estimation asymptotique du déterminant de $T_{\Lam}(f)$.
 \end{abstract}
 \selectlanguage{english}
 \begin{abstract}
We consider the class of positive bounded and semi-continuous functions defined on the torus $\tore^2$. 
If $f$ belongs to this class, $f$  will be considered as the symbol of a Toeplitz operator  truncated on a triangle $\Lam$ parametrised by an integer number $\lambda$. This operator is denoted by $T_{\Lam}(f)$. We develop a geometric structure of the inverse of $T_{\Lam}(f)$ and give an asymptotical development of the trace of $T_{\Lam}(f)^{-1}$ wich brings out the geometry of the triangle. The foundation of this result consists in the possibility of $f$ having  a factorisation of type $\alpha\overline \alpha$ where the spectrum of $\alpha$ will be localised in a given semi-cone. This trace theorem allows in particular to find again the Linnik-Szegö theorem about the asymptotical evaluation of the determinant of $T_{\Lam}(f)$.
 \end{abstract}
 \selectlanguage{francais}
\section{Introduction}
Cet article s'inscrit dans une longue tradition de travaux prolongeant le théorème limite de Szegö, tel qu'il appara\^{i}t dans \cite{SZE} (page $194$) ou dans \cite{GS} (page $65$). Ce travail est  plus  particulièrement dans la continuité du second théorème de Szegö  tel qu'on le trouve dans \cite{SZO} où appara\^{i}t un opérateur de Toeplitz tronqué sur un intervalle de $\rr$ défini par un symbole qui est une fonction définie sur le tore $\tore^1$.  Ce théorème a été généralisé en dimension $1$ par Kac, Murdock et Szegö dans \cite{KMS}.

 En dimension $1$, les extensions du théorème de Szegö passent par le choix de symboles singuliers.  Par exemple citons dans ces travaux \cite{RaSegh} où l'on se donne un symbole factorisé sous la forme
$ f=|1-e^{i\theta}|^{2\alpha} f_{1},\;\alpha\not= 0,\alpha\in]-1/2,1/2[$ et $f_{1}$ régulière, mais également \cite{RamRin}, où 
$f=|1-e^{i\theta}|^2 f_{1}$, où $f_{1}^{\pm1}$ est analytique au voisinage de $0$. Le cas $\alpha=1$ est aussi l'objet d'articles tels que \cite{Part} et \cite{JMR01}. Ces factorisations mettent en évidence l'ordre du zéro du symbole et établissent des propriétés qui sont liées à l'ordre de ces zéros.

Cependant,  certaines applications imposent un symbole spécifique et c'est alors par des théorèmes que sont décrites ses factorisations. Dans ce cadre et toujours en dimension $1$, citons le cas où $f=1-\Phi$, la fonction $\Phi$ étant la fonction génératrice d'une variable aléatoire $X$. Les probabilités imposent ici des hypothèses sur $X$ et ces hypothèses donnent lieu, éventuellement, à une factorisation du symbole. On peut consulter sur cet aspect les articles \cite{De2Ri2} et \cite{De3Ri3}. En outre 
on montre dans \cite{RaRi}, l'importance d'un théorème de trace pour une marche aléatoire. On pourra consulter le livre de F. Spitzer (\cite{FS}) pour le lien entre les marches aléatoires et les matrices de Toeplitz.

Revenons maintenant au problème en dimension supérieure à un. Le passage de la dimension un à la dimension deux se fait naturellement avec des techniques identiques à celles de la dimension un, si l'on envisage des symboles matriciels à coefficients de Fourier matriciels. On suppose, dans ce cadre, une factorisation du symbole matriciel $F$ sous la forme $F=GG^*$ avec des hypothèses spécifiques sur $G$. On trouvera cette approche dans \cite{JC2}. Cependant cette approche ne peut constituer une méthode générale d'une théorie multidimensionnelle car elle ne peut traiter que le cas d'une sous-classe d'opérateurs tronqués sur un (multi-)rectangle.

\noindent L'approche multidimensionnelle, permettant une extension à des domaines convexes autres que les rectangles, concerne pour l'essentiel des symboles réguliers.
Un des théorèmes principaux de ce travail donne un développement asymptotique de la trace de l'inverse d'un opérateur de Toeplitz tronqué sur un triangle, opérateur défini par un symbole positif.
 L'origine des idées de ce travail réside dans deux articles de Widom et Linnik.  Ce dernier, dans  \cite{LINN}, démontre le théorème limite de Szegö dans le cadre multidimensionnel discret. Sa méthode s'appuie sur l'estimation des puissances de matrices de Toeplitz. Widom, quant à lui, établit dans \cite{Wiwi} à partir du logarithme du déterminant de la matrice de Toeplitz un analogue du théorème de Szegö, en dimension $n$ dans le cas d'un opérateur  continu (opérateur de convolution), en s'appuyant sur une approximation fine de l'inverse de la matrice de Toeplitz. Suite à ces deux approches, l'un des auteurs de ce travail  a proposé une inversion exacte dans le cadre multidimensionel pour un opérateur de Toepltiz tronqué sur un multirectangle et en a déduit un théorème de trace de l'inverse de cet opérateur avec un symbole de la forme $1/|P|^2$, où $P$ désigne un polynôme trigonométrique (voir  \cite{sghr}). Ce théorème de trace  a été généralisé par B.H. Thorsen dans \cite{THO}, qui prend pour symbole $f=gh$ où $g^{\pm1}$ et $\bar h^{\pm1}$ sont dans un espace de Hardy $H^\infty(U^n)$, l'opérateur étant toujours tronqué sur un 
multi-rectangle. Dans \cite{SeKa}, la factorisation étant  $1/|P|^2$, avec $P$ polynôme trigonométrique, les auteurs construisent une démarche formelle pour des opérateurs tronqués sur des polytopes.

   Nous allons donner, dans cet article une construction précise de l'inverse pour aboutir à une formule de trace. Le résultat est obtenu à partir de l'existence d'une factorisation \og naturelle \fg\; du symbole associée à la géométrie du polygone convexe sur lequel l'opérateur est tronqué. Les hypothèses exigées concernant le symbole $f$ sont à peine plus fortes que celles requises pour qu'une fonction positive sur $\tore^2$ soit le module d'une limite radiale d'une fonction de l'espace de Hardy $H^\infty(U^2)$ (voir \cite{RUD} page $55$.)
   Il y a plusieurs factorisations possibles du symbole pour aboutir à une formule d'inversion d'un opérateur de Toeplitz tronqué sur un polygone convexe. Le choix d'une factorisation influence la taille des opérateurs de Hankel qui interviennent dans le calcul de l'inverse. Dans cet esprit, nous montrons qu'il existe une factorisation optimale appelée ici \og factorisation  minimale\fg. 
  Le reste de l'article est consacré au théorème de trace de l'inverse de l'opérateur de Toeplitz tronqué sur un triangle,  à son corollaire, le théorème de Linnik (\cite{LINN}), avec des hypothèses un peu différentes sur la régularité du domaine de troncature (Linnik exige un bord $C^2$ par exemple et des hypothèses de régularité plus fortes sur le symbole) qui donne une évaluation asymptotique du déterminant de ce m\^{e}me opérateur. 
  Notons que la formule de trace a ceci de remarquable, c'est qu'elle s'exprime en fonction de mesures liées au symbole et à la géométrie du triangle. Le lien entre la géométrie du polygone et les coefficients de la trace est plus apparent que pour le rectangle dont les angles aux sommets sont fixés.
  
 Indiquons  pour finir que la méthode employée pour le triangle est généralisable à une polygone convexe, objet d'un 
  autre travail.
  La section suivante contient les notations utilisées dans ce papier, bien que la plupart d'entre elles  soient standards. En section \ref{resu}, nous détaillons les principaux résultats évoqués dans l'introduction.

\section{Notations}\label{nonotes}
\emph{On trouvera  les notations suivantes dans la section \ref{cocou} et \ref{lap}.}
\begin{enumerate}
\item $U$ est le disque unitaire de $\cc$ et $\tore$ est le bord de $U$.
\item $m_{2}$ est la mesure de Lebesgue sur $\tore^2$.
\item  $H(U^2)$ désigne les fonctions holomorphes sur $U^2$ et
$H^\infty(U^2)$ désigne les fonctions holomorphes bornées sur $U^2$.
\item Si $\mu$ est une mesure complexe sur $\tore^2$, $P[d\mu]$ désigne son intégrale de Poisson sur $U^2$. De fa\c con explicite, posant
$z_{j}=r_{j}e^{i\theta_{j}}, |r_{j}|<1$, $w_{j}=e^{i\varphi_{j}}$, $P_{r_{j}}(\theta_{j}-\varphi_{j})=\frac{1-r_{j}^2}{1-2r_{j}\cos(\theta_{j}-\varphi_{j})
+r_{j}^2}$,

\noindent $P(z,w)=\prod_{j=1}^2 P_{r_{j}}(\theta_{j}-\varphi_{j}),z=(z_{1},z_{2}), w=(w_{1},w_{2})$, on a $P[d\mu](z)=\int_{\tore^2}P(z,w)d\mu(w)$
\item $Y_{2}$ désigne le cône $\zz_{+}^2\cup (-\zz_{+}^2)$ et $\tilde{Y}_{2}$ son image par la symétrie de $\zz^2$: $(m,n)\mapsto (-m,n)$.
\item  Pour toute fonction $g\in L^1(\tore^2)$, on pose $\hat{g}=\{\hat{g}(k)\}_{k\in \zz^2}$ la famille des coefficients de Fourier de $g$  et $\spec({g})=\{k\in\zz^2;\;\hat{g}(k)\not=0\}$.
\item Si $C$ est un cône de $\rr^2$ de sommet $O=(0,0)$, inclus dans $Y_{2}$, on pose $$C_{+}=C\cap\zz_{+}^2\;\text{et}\;C_{-}=C\cap(-\zz_{+}^2).$$
Si $C$ est le cône $\{\alpha e_1+\beta e_2; (\alpha,\beta)\in\rr_+^2\cup(-\rr_+^2)\}$, on le note $C [e_1, e_2]$. 
$C^+[e_1, e_2]=\{\alpha e_1+\beta e_2; (\alpha,\beta)\in\rr_+^2\}$ sera aussi noté $C^+[d_1^+, d_2^+]$ si $ d_{i}^+=\{\alpha e_{i}\;;\;\alpha\ge 0\}$
\item Si $\alpha=(\alpha_{1},\alpha_{2}) \in \rr^2$ et $\theta=(\theta_{1},\theta_{2})\in \rr^2$, on pose $ \alpha.\theta=\sum_{i=1}^2\alpha_{i}\theta_{i}$, ainsi que $ e^{i\theta}=(e^{i\theta_{1}},e^{i\theta_{2}})$. 
\item Si $u$ est une fonction définie sur $U^2$, on pose avec la notation du point précédent $u^*(e^{i\theta})=\lim_{r\to 1}u(re^{i\theta})$.
\end{enumerate}

\noindent\emph{ On trouvera  les notations suivantes dans les sections \ref{resu}, \ref{polo} et \ref{lap}}
\begin{enumerate}
\item Pour tout sous-ensemble $A$ de $\zz^2$, $\Po(A)$ est  l'espace vectoriel engendré par les  
polyn\^omes trigonom\'etriques $g$ tels que $\spec({g})\subset A$ et 
$\Pi_A$ le projecteur orthogonal de $L^2(\tore^2)$ dans
$\Po(A)$.
\item On note $m_{2}$ la mesure de Lebesgue sur le tore $\tore^2$ et $\la,\ra$ désigne le produit scalaire de $L^2(\tore^2)$. Si $f$ est une fonction mesurable positive et bornée sur le tore, $ L_{1/f}^2(\tore^2)$ désigne l'espace des fonctions complexes de carré intégrable sur $\tore^2$, muni du produit scalaire $\langle g,h\rangle=\int_{\tore^2}g\bar{h}d\sigma$ où $d\sigma=\dfrac{1}{f}dm_{2}$.
\item le signe $\equiv$ désigne une égalité par définition. \item Si $E$ est un ensemble fini, $|E|$ désigne son cardinal.

\end{enumerate}

\section{Les principaux résultats}\label{resu}
\subsection{ Les principales notations et hypothèses liées au triangle et au symbole.}\label{trisym}
\begin{enumerate}
\item \emph{Hypothèses liées au triangle.}

 \'{E}tant donné un paramètre entier et positif, $\lambda$, on considère dans cette section le triangle $\Lam_{\lambda}$,  formé par des points $O,A_{1}A_{2}$ à coordonnées dans $\nn^2$. On note $c_{1},c_{2},c_{3}$ les côtés $OA_{1}$, $A_{1}A_{2}$, $A_{2}O$, $\nu_{i}=\{(\alpha_{i},\beta_{i})\in\zz^2\}_{i=1,2}$ des vecteurs normaux aux trois côtés $c_{i}$, à composantes entières et \emph{premières entre elles}.  On pose $O=(0,0),A_{1}=(\lambda\beta_{1},-\lambda\alpha_{1}), A_{2}=(a\lambda,0)$ où $a$   entier naturel fixé et pour finir  $\nu_{3}=(0,-1)$ (voir figure \ref{triangle}).
 Notons $l_{i}(\lambda)$ la longueur (euclidienne) du côté $c_{i}$  du triangle $\Lam_{\lambda}$ de la figure \ref{triangle}, et $n_{i}=(a_{i},b_{i})=\nu_{i}/||\nu_{i}||$ le vecteur unité normal au côté $c_{i}$ et posons 
$\mathfrak{S}_{1}(\lambda)=\frac{1}{2}\sum_{i=1}^3(-1)^{i}a_{i}l_{i}(\lambda)$, $\mathfrak{S}_{2}(\lambda)=\frac{1}{2}\sum_{i=1}^3(-1)^{i}b_{i}l_{i}(\lambda)$, $\mathfrak{S}_{1}=\mathfrak{S}_{1}(1),\mathfrak{S}_{2}=\mathfrak{S}_{2}(1)$
  \begin{figure}
\unitlength 0.9mm
\begin{picture}(105, 85)
\thinlines
 \put(20,0){\line(0,1){85}}
 \put(0,20){\line(1,0){105}}
 \put(0,0){\line(1,1){85}}
  \put(80,20){\vector(0,-1){15}}
 \put(50,50){\vector(-1,1){10}}
 \put(80,50){\vector(2,1){10}}
\thicklines
  \put(20,20){\line(1,0){75}}
  \put(20,20){\line(1,1){50}}
  \put(70,70){\line(1,-2){25}} 
 \put(92,55){$\vec\nu_2$}
 \put(81, 5){$\vec\nu_3$}
 \put(42, 61){$\vec\nu_1$} 
 \put(17,21){$0$}
 \put(93,17){$A_2$}
 \put(66,72){$A_1$}
 \put(45,17){$c_3$}
 \put(40,38){$c_1$}
 \put(70,60){$c_2$}
 \put(5,12){$\mathcal C^-$}
 \put(30,25){$\mathcal C^+$}
 \end{picture}
 \caption{triangle}\label{triangle}
 \end{figure}
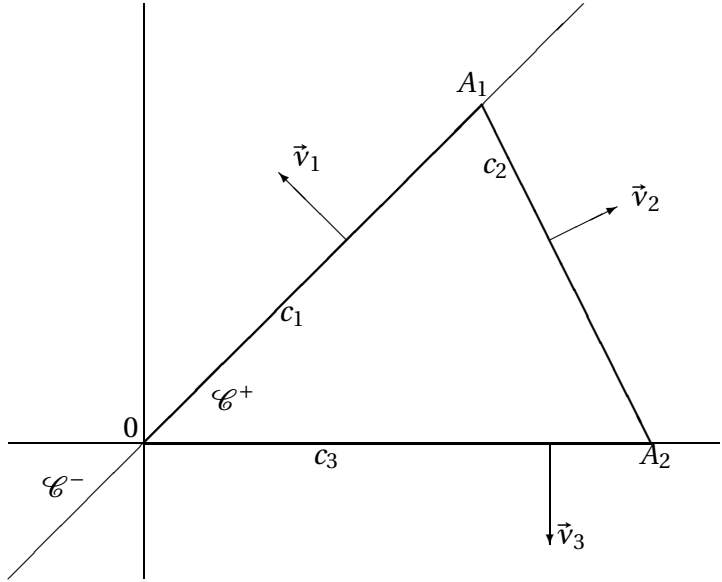
Supposer que l'angle du sommet $A_{2}$ de la figure \ref{triangle} est aigu, 
l'angle du sommet $O$ ayant une mesure en radian dans $[0,\Pi/2]$, équivaut à l'hypothèse $(\mathfrak{T})$ définie par  les inégalités 
\begin{equation}\label{posi}
\mathfrak{S}_{1}>0,\mathfrak{S}_{2}>0
\end{equation}
Remarquons que l'hypothèse $(\mathfrak{T})$  ne restreint pas la généralité du triangle à un déplacement près. 
Avec l'hypothèse (\ref{posi}), on considère l'espace vectoriel $L^\infty(\tore^2)/\sim$ où $\sim$ est la relation d'équivalence sur $L^\infty(\tore^2)$ définie par $f\sim g$ si $f-g$ est une constante. Nous noterons  $\ma L^\infty(\tore^2)$ l'espace $L^\infty(\tore^2)/\sim$. On considère le sous-espace vectoriel $\mathfrak K$ de   $L^\infty(\tore^2)$ défini par 
$\mathfrak K=\{g\in \ma L^\infty(\tore^2)\;;\;||g||^2_{\mathfrak K,\Lam}=\sum_{(u,v)\in\zz^2}(\mathfrak{S}_{1}|u|+\mathfrak{S}_{2}|v|)|\widehat g (u,v)|^2<\infty$. 
Alors $(\mathfrak K,||.||_{\mathfrak K,\Lam})$ est un espace vectoriel normé. Sa norme est déterminée par le produit scalaire :
\begin{center}
 $\la g_{1},g_{2}\ra=\sum_{(u,v)\in\zz^2}(\mathfrak{S}_{1}|u|+\mathfrak{S}_{2}|v|)\widehat g_{1} (u,v)\overline{\widehat g_{2} (u,v)}$
\end{center}
\item \emph{Hypothèses liées au symbole.}

Soit $f$ le symbole. C'est  une fonction sur $\tore^2$ sur laquelle nous aurons à faire une  des trois  hypothèse suivantes. 
\begin{enumerate}
\item Hypothèse $(\ma{H}_{1})$ :

$f$ une fonction positive, essentiellement bornée et semi-continue inférieurement sur $ \tore^2$. 

\item Hypothèse $(\ma{H}_{2})$ :

$f$ vérifie $(\ma{H}_{1})$ et $1/f$ est essentiellement bornée.

\item Hypothèse $(\ma{H}_{3})$ :

$f$ vérifie $(\ma{H}_{2})$ et les fonctions $f$ et $\ln f$ appartiennent à $\mathfrak{K}$.
\end{enumerate}
\end{enumerate}
\subsection{Opérateur de Toeplitz tronqué sur un polygone convexe, de symbole régulier}
\begin{definition}\label{eccolo}
Soit $A$ un polygone convexe de $\zz^2$. Si $f$ est une fonction positive de $L^\infty(\tore^2)$ et $T(f)$ la multiplication
par $f$ dans $L^2(\tore^2)$, alors l\emph{'op\'erateur de Toeplitz tronqu\'e sur} $\Po(A)$ associ\'e au symbole $f$ et not\'e $T_A(f)$est d\'efini par 
$$\forall q\in \Po(A),\; \;T_A(f)(q)=\Pi_A (fq).$$
 \end{definition}

 Ce théorème décrit le comportement asymptotique du déterminant de l'opérateur de Toeplitz tronqué sur le triangle $A=\Lam_{\lambda}$ lorsque $\lambda$ tend vers $+\infty$.
 
 \subsection{Théorème d'inversion et théorème de  trace. Cas du triangle.}
 Sous les hypothèses $(\ma{H}_{1})$, la matrice de Toeplitz 
 $T_{\Lam_{\lambda}}(f)$ est inversible et son inverse est donné pour tout polynôme $q$ de $\ma{P}(\Lam_{\lambda})$ par une égalité de la forme
 $T_{\Lam_{\lambda}}(f)^{-1}(q)=\frac{q}{f}+\mathfrak{H}(q)$ où $\mathfrak{H}$ est un opérateur lié à la description de $\Lam_{\lambda}$ comme intersection de demi-espaces. 
  \begin{theorem}{ théorème de trace}\label{otlitchno}
 
On suppose que $f$ vérifie l'hypothèses $(\ma{H}_{2})$. Alors 
\begin{align*}
\trace\left(T_{\Lam_{\lambda}}(f)^{-1}\right)=|\Lam| \;||\frac{1}{f}||_{1}&+\left(\sum_{(u,v)\in\zz^2}|u|\;\widehat{\log \left(\frac{1}{f}\right)}(u,v)\overline{\widehat{\frac{1}{f}}(u,v)}\right)\mathfrak{S}_{1}(\lambda)\\&+\left(\sum_{(u,v)\in\zz^2}|v|\;\widehat{\log \left(\frac{1}{f}\right)}(u,v)\overline{\widehat{\frac{1}{f}}(u,v)}\right)\mathfrak{S}_{2}(\lambda)+o(\lambda).
\end{align*}
\end{theorem}
On en déduit le comportement asymptotique du déterminant de l'opérateur de Toeplitz tronqué sur le triangle $\Lam_{\lambda}$ lorsque $\lambda$ tend vers $+\infty$ en retrouvant dans le cadre du triangle le théorème de 
Linnik.
\begin{coro}{théorème du déterminant de Szegö-Linnik}\label{determi}
 
Soit $\Lam_{\lambda}$ le triangle décrit en section \ref{trisym}. Faisons l'hypothèse $(\mathfrak{T})$ sur le triangle et l'hypothèse $(\ma{H}_{3})$ sur le symbole. 

\noindent Posons 
$\mu_{1}(f)=-\frac{1}{2}\sum_{(u,v)\in\zz^2}|u| |\widehat{\ln f}(u,v)|^2$,
$  \mu_{2}(f)=-\frac{1}{2}\sum_{(u,v)\in\zz^2}|v| |\widehat{\ln f}(u,v)|^2$.

 Alors
 
  \[\det\;T_{\Lam_{\lambda}}(f)=e^{|\Lam_{\lambda}|\;||\ln f||_{1}}e^{-\lambda\big(\mathfrak S_{1}\mu_{1}(f)+\mathfrak S_{2}\mu_{2}(f)+o(1)\big)}.\]
\end{coro}

  \subsection{Théorème de factorisation}
  Ce théorème montre que sous des hypothèses raisonnables de régularité, une fonction positive $f$ sur le tore $\tore^2$ se factorise en un produit de deux fonctions conjuguées à support dans un cône donné. Ces hypothèses sont contenues dans les hypothèses des théorèmes de trace et du déterminant.
  Plus précisément :
  \begin{theorem}\label{clef}
Soit $f$ une fonction positive, essentiellement bornée et semi-continue inférieurement sur $\tore^2$ et soit  $\C$ un cône de sommet $O$ inclus dans $Y_{2}$. Il existe une fonction  $\alpha\in H^{++}=\{\psi\in L^\infty(\tore^2);\;\forall k\notin \nn^2,\hat{\psi}(k)=0\}$ telle que $f=\alpha\bar{\alpha},\;\spec({\alpha})\subset \C_{+}$ et $\spec({\alpha^{-1}})\subset \C_{+}=\C\cap\zz_{+}^2$. 
\end{theorem}
\begin{coro}\label{clebar}
Soit $f$ une fonction positive, essentiellement bornée et semi-continue inférieurement sur $\tore^2$ et soit  $\C$ un cône de sommet $O$ dont on note $\C^+$ un des demi-cônes. On pose $H^{\C^+}=\{\psi\in L^\infty(\tore^2);\;\forall k\notin \C^+,\hat{\psi}(k)=0\}$.   Il existe une fonction  $\alpha\in H^{\C^+}$ telle que $f=\alpha\bar{\alpha}$.  
\end{coro}
Le corollaire \ref{clebar} n'est pas utilisé dans la démonstration du théorème de trace qui applique uniquement le théorème \ref{clef}. Cependant il est utile dans le théorème d'inversion sous sa forme générale.
Indiquons que le théorème \ref{clef} généralise un résultat de Rudin (voir\cite{RUD}) qui dit que si $f$ est semi-continue inférieurement sur $\tore^2$ et si $f$ est dans $L^1(\tore^2)$ alors $P[f-d\sigma]$ est une fonction à partie réelle holomorphe pour une certaine mesure positive singulière  sur $\tore^2$.
\subsection{Plan des démonstrations}
On commence par établir le théorème de factorisation en section \ref{cocou}. Puis on introduit la notion de factorisation relative à un polygone convexe dans la section \ref{polo}. On établit alors un théorème d'inversion de 
l'opérateur de Toeplitz tronqué sur un polygone convexe (théorème \ref{invon}) sous les hypothèses du théorème \ref{clef}. 
Le théorème d'inversion est à la base du théorème de trace. Celui-ci est démontré en section \ref{lap} lorsque le polygone est un triangle. Le corollaire \ref{determi} est démontré en section \ref{deter}. 
\section{Factorisation relative à un cône\label{cocou}}
La démonstration du théorème \ref{clef} passe par les lemmes suivants.
\begin{lem}\label{cone}
Soit $\C$ un cône, non réduit à une droite, de sommet $O$ et $S$ une partie finie de $\zz^2$. Il existe $v(\C,S)$ dans $ \C$ tel que 
\begin{enumerate}
\item $\forall n\in\zz^*,\;O\notin S-n v(\C,S)$,
\item $\forall n\in\zz^*,\; S-nv(\C,S)\subset \C$.
\end{enumerate}
\end{lem}
\begin{preuve}{du lemme \ref{cone}}
	
\noindent Preuve en deux parties.

 \emph{Première partie} : \emph{on suppose que le cône $\C$ est engendré par une $\zz$-base $(e_{1},e_{2})$ de $\zz^2$ c'est à dire que 
$\C=\C_{+}\cup (-\C_{+})$ avec $\C_{+}=\{n e_{1}+m e_{2};(n,m)\in\nn^2\}$.} 

Dans ce cas, supposons tout d'abord que $S$ est réduit à un point.
Si $S=\{O\}$, c'est évident. Sinon tout  singleton $S$,  pouvant s'écrire sous la forme $S=\{ae_{1}+be_{2}\}$, montrons qu'il existe $(\lambda,\mu)\in Y_2\setminus \{O\}$ tel que $a e_1+b e_2-n(\lambda e_1+\mu e_2)\in Y_2\setminus \{(0,0)\}$, pour tout $n\in Z^*$.
Si $(a,b)\not=(0,0)$, on prend $(\lambda,\mu)=(2|a|,2|b|)$ et
si $(a,b)=(0,b)$ où $b\not= 0$, on prend $(\lambda,\mu)=(|b|/b,-2b)$. 
Alors dans les deux cas, pour tout $n\in\zz^*$, $
	 (a-n\lambda)(b-n\mu)>0 
$. 

\noindent Supposons maintenant que $S$ soit fini de cardinal $n$.
On écrit $S=\{a_i\}_{i\in\{1,\ldots,n\}}$. Pour tout $i\in\{1,\ldots,n\}$, 
il existe $v_i\in\C[e_1,e_2]$ pour lequel les deux propriétés du lemme \ref{cone} sont vérifiées. Alors $v(\C[e_1,e_2],S)=\sum_{i=1}^n v_{i}$ répond à la question.

\emph{Deuxième partie : le cône $\C$ est quelconque}.

Alors il existe un cône $\C[e_{1},e_{2}]$ engendré par une $\zz$-base de $\zz^2$ tel que $\C[e_{1},e_{2}]\subset\C$. Supposons en effet que 
$\C=\C[(a,b),(c,d)]$. On ne restreint pas la généralité du résultat en supposant $a,b,c,d$ positifs et alors $\C_{+}$ est inclus strictement dans 
$\nn^2$. La densité des rationnels fournit un couple $(m,n)$ d'entiers naturels tel que $\frac{b}{a}<\frac{n}{m}<\frac{d}{c}$. Il existe un couple d'entiers positifs $(u_{0},v_{0})$, premiers entre eux tels que $mu_{0}-nv_{0}=\pm1$. Mais alors l'ensemble des couples $(u,v)$ vérifiant 
$mu_{}-nv_{}=\pm1$ est de la forme $u(k)=u_{0}+kn$ et $v(k)=v_{0}+km$ où $k\in\zz$. D'où $\frac{v}{u}=\frac{n}{m}+o(1))$ avec $\lim_{k\infty}o(1)=0$.  Ceci montre que 
$(u(k),v(k))\in\C_{+}$ lorsque $k$ est assez grand et on prend dans ce cas $e_{1}=(m,n), e_{2}=(u(k),v(k))$.
\end{preuve}
\begin{lem}\label{sing}
Soit $f\in L^1(\tore^2)$ une fonction positive semi continue inférieurement et $\C$ un cône de sommet $O$ inclus dans $Y_{2}$. Il existe une fonction holomorphe $g$ et une mesure positive $\sigma$ singulière sur $\tore^2$ (c'est à dire $\sigma\bot m_{2}$) vérifiant :
\begin{enumerate}
\item $P[f-d\sigma]=\re(g)$,
\item $spec ({g})\subset \C_{+}.$
\end{enumerate}
\end{lem}
\begin{preuve}{ en deux étapes}

\noindent\textbf{Première étape} : \emph{montrons que si $f$ est un polynôme trigonom\'etrique positif ou nul, il existe une mesure singuli\`ere sur $\tore^2$ telle que pour tout $a\notin \C_{+}\cup (-\C_{+})=\C\cap Y_{2}$ on a $ \hat{\sigma}(a)=\hat{f}(a)$.}

\noindent Puisque $f$ est un polynôme trigonom\'etrique, $\hat{f}(k)=0$ en dehors d'un ensemble fini $S$. On peut donc trouver, d'après le lemme \ref{cone} un vecteur $v(\C,S)$ noté plus plus simplement ici $v=(v_{1},v_{2})$ tel que : 
\begin{equation}\label{ere}
\left\{
\begin{array}{l}
\forall  p\in\zz^*,\;0\notin S-pv,\\
\forall  p\in\zz^*,\;S-pv\subset \C \;\cap Y_{2}.
\end{array}
\right.
\end{equation}
Soit $H$ le sous-groupe de $\tore^2$ d\'efini par $H=\{e^{i\theta}=(e^{i\theta_{1}},e^{i\theta_{2}})\in \tore^2; e^{i\theta.v}=1\}$.
On note $\mu_{H}$ la mesure de Haar associ\'ee au groupe $H$. Comme $m_{2}(H)=0$, on a $\mu_{H}\perp m_{2}$. On consid\`ere alors la mesure singuli\`ere $d\sigma=fd\mu_{H}$. Si on note $|\sigma|$ la variation totale de $\sigma$ et $||\sigma||=\int_{\tore^2}d|\sigma|$ sa norme, on a $||\sigma||=||f||_{1}=\int_{\tore^2}fd\mu_{H}$, car $f\ge 0$ (voir \cite{RU}). \'{E}valuons les coefficients de Fourier de $\sigma$. On a 
$\hat{\sigma}(a)=\int_{\tore^2}e^{-i\theta.a}\left(\sum_{finie}\hat{f}(k)e^{ik.\theta}\right)d\mu_{H}(\theta)=\sum_{finie}\hat{f}(k)\int_{\tore^2}e^{i\theta.(k-a)}d\mu_{H}(\theta)$. Or l'intégrale $\int_{H}e^{i\theta.(k-a)}d\mu_{H}(\theta)$ est nulle si $k-a$ n'est pas un multiple entier de $v$. En effet, $\mu_{H}$ est invariante par translation. Donc si $e^{i\alpha}\in H$ et si $k-a$ n'est pas un multiple de $v$, alors
$\int_{\tore^n}e^{i\theta.(k-a)}d\mu_{H}(\theta)=\int_{H}\underbrace{e^{i\theta.(k-a)}}_{{\varphi(\theta)}}d\mu_{H}(\theta)=\int_{H}\varphi(\theta+\alpha)d\mu_{H}
(\theta)=\int_{H}e^{i(\theta+\alpha).(k-a)}d\mu_{H}(\theta)=e^{i\alpha(k-a)}\int_{H}e^{i\theta.(k-a)}d\mu_{H}(\theta).$ Or $e^{i\alpha.(k-a)}\not=1$, ce qui permet 
de conclure. Par contre si $k-a$ est un multiple entier de $v$, alors $e^{i\theta.(k-a)}=1$. Ainsi,
$$\hat{\sigma}(a)=\sum_{finie}\hat{f}(k)\left( \underbrace{\int_{H}e^{i\theta.(k-a)}d\mu_{H}(\theta)}_{=1}+\underbrace{\int_{\tore^n\setminus H}e^{i\theta.(k-a)}d\mu_{H}(\theta)}_{=0}\right)$$
car $\mu_{H}(H)=1$ et le support de $\mu_{H}$ est inclus dans $H$. On en d\'eduit que $\hat{\sigma}(a)=\sum_{j=-\infty}^{+\infty}\hat{f}(a+jv)=
\hat{f}(a)$ si $a\notin \C_{+}\cup (-\C_{+})$, d'après la condition (\ref{ere}). 

\noindent\textbf{Deuxi\`eme \'etape} : \emph{montrons que si $f\in L^1(\tore^2)$ est une fonction positive ou nulle, semi-continue inf\'erieurement il existe une mesure singuli\`ere sur $\tore^2$ telle que pour tout $a\notin \C_{+}\cup (-\C_{+})$ on a $ \hat{\sigma}(a)=\hat{f}(a)$.}

La fonction $f$ étant semi-continue inférieurement,  on peut \'ecrire $f=\sum_{j=0}^\infty f_{j}$ o\`u $f_{j}$ est positive ou nulle et continue (voir \cite{SCHW}). Et on a \'egalement l'\'egalit\'e 
$f_{j}=\sum_{i=0}^\infty g_{ij}$ o\`u les $g_{ij}$ sont des polynômes trigonom\'etriques positifs ou nuls. En regroupant ce qui pr\'ec\`ede on peut donc \'ecrire $f$ sous la forme $f=\sum_{j=0}^\infty g_{j}$ avec $g_{j}$ polynôme trigonom\'etrique positif ou nul. D'apr\`es la premi\`ere \'etape, on peut associer \`a chaque $g_{j}$ une mesure singuli\`ere $\sigma_{j}$ v\'erifiant :
\begin{itemize}
\item  $||\sigma_{j}||=||g_{j}||_{1}$,
\item  $\forall  a\notin  \C_{+}\cup (-\C_{+}),\;\hat{\sigma}_{j}(a)=\hat{g_{j}}(a)$.
\end{itemize}
De plus $\sum_{i}||g_{i}||_{1}=\sum_{i}\int_{\tore^2}|g_{i}|d\mu_{H}=\int_{\tore^2}(\sum_{i}g_{i})d\mu_{H}=\int_{\tore^2}fd\mu_{H}$ car $g_{i}\ge 0$. Ainsi la s\'erie $\sum_{i}\sigma_{i}$ converge normalement, donc converge vers la mesure complexe $\sigma=\sum_{i}\sigma_{i}$ : en effet l'espace vectoriel normé  des mesures complexes est complet puisque  isométrique au dual topologique de $C(\tore^2)$ d'après un théorème de représentation de Riesz. La mesure $\sigma$ est une mesure complexe singuli\`ere. Si $k\notin \C_{+}\cup (-\C_{+})$, on a $\hat{\sigma}(k)=\sum\hat{\sigma}_{i}(k)=\hat{f}(k)$.

La démonstration s'achève avec la remarque suivante : Si $\mu$ est une mesure complexe sur $\tore^2$, $h=P[d\mu]$ est la partie réelle d'une fonction holomorphe $g$ si et seulement si $\hat{\mu}(k)=0$ si $k\notin \nn^2\cup (-\nn^2)$ et alors $\spec({h})=\spec({\mu})$ (voir \cite{RUD}, théorème 2-4-1). Or si $\spec {(\re(g)})\subset \C_{+}\cup (-\C_{+})$, alors  $g$ étant holomorphe le support de $\hat{g}$ est nécessairement inclus dans $\C_{+}$. Ainsi $P[f-d\sigma]$ est elle la partie réelle 
d'une fonction holomorphe à support dans $\ma C_{+}$.
\end{preuve}
\begin{lem}\label{rude}
Notons $H^{++}=\{u\in L^\infty(\tore^2); \forall k\notin\nn^2,\hat{u}(k)=0\}$. Si $u\in H^{++}$ vérifie l'inclusion $\spec({u})\subset \C_{+}$ alors pour tout $n\in \nn$, on a : $\spec({u^n})\subset \C_{+}$. 
\end{lem}
\begin{coro}[du lemme \ref{rude}]\label{raide}
Soit $u$ une fonction de $ H^{++}$ telle que  $\spec({u})\subset \C_{+}$. Alors $\spec({\exp(u)})\subset\C_{+}$.
\end{coro}

\begin{preuve}{du lemme \ref{rude} et du corollaire \ref{raide}}
Le résultat étant vrai pour $n=1$, il suffit pour mettre en \oe uvre le principe de récurrence de montrer que si $u,v\in H^{++}$ vérifient $\spec({u})\subset \C_{+}$ et $\spec({v})\subset \C_{+}$, alors  $\spec({uv})\subset \C_{+}$. On se rappelle qu'il y a une isométrie bijective $\Phi$ entre $H^\infty(U^2)$ et $H^{++}$ définie pour tout $\alpha\in\ H^\infty(U^2)$ par $\Phi(\alpha)=\alpha^*$ : ceci est une  de ces propriétés sur l'existence de limites radiales de l'analyse des fonctions holomorphes  à une variable qui s'étendent naturellement  aux fonctions holomorphes à plusieurs variables (voir \cite{RU} et \cite{RUD}). Il existe donc deux fonctions $\varphi$ et $\psi$ dans $H^\infty(U^2)$ telles que pour presque tout $\theta$ au sens de la mesure de Lebesgue sur $\tore^2$, on ait : $u=\varphi^*$ et $v=\psi^*$.
On a pour $\varphi$ et $\psi$ respectivement les développements en séries entières 
$\varphi=\sum_{(i,j)\in\nn^2}a_{ij}z_{1}^{i}z_{2}^j$ et $\psi=\sum_{(i,j)\in\nn^2}b_{ij}z_{1}^{i}z_{2}^j$  qui convergent normalement sur tout sous-polydisque de $U^2$ et par conséquent $$\varphi\psi=\sum_{(s,t)\in\nn^2} (\sum_{\genfrac{}{}{0pt}{1}{i+k=s} {j+l=t}}{a_{ij}b_{kl}})z_{1}^sz_{2}^t=\sum_{(s,t)\in\nn^2} c_{st}z_{1}^sz_{2}^t.$$
De plus pour tout $(i,j)\in\C_{+}$ et tout $(k,l)\in\C_{+}$, on a $(s,t)=(i,j)+(k,l)\in\C_{+}$. Le théorème d'Abel permet d'identifier les coefficients 
$c_{st}$ comme les coefficients de Fourier de $uv=\varphi^*\psi^*$, ce qui prouve que $\spec({uv})\subset \C_{+}$. Et on conclut.
Le corollaire s'en déduit par le développement en série entière de l'exponentielle.
\end{preuve}
\begin{lem}\label{unsur}
Soit $u$ une fonction non nulle de $ H^{++}$ telle que  $\spec({u})\subset \C_{+}$. Alors $\spec({u^{-1}})\subset\C_{+}$.
\end{lem}
\begin{preuve}{du lemme \ref{unsur}}
L'argumentation est du m\^{e}me type que celle du lemme \ref{rude}. Si $u\not=0$ et $u\in H^{++}$, alors il existe $\alpha\in H^\infty(U^2)$ 
telle que $u=\alpha^*$ et $\alpha\not= 0$ sur une couronne $C_{\varepsilon}=\{(z_{1},z_{2});1-\varepsilon<|z_{i}|\le 1,i=1,2\}$ pour $\varepsilon$ assez petit. Alors sur $\tore^2$, on a  $u^{-1}={\alpha^*}^{-1}={\alpha^{-1}}^*$. On peut associer à $\alpha$ une série formelle dont on peut supposer le terme constant égal à $1$, de la forme $S(X,Y)=1+T(X,Y)$ avec $T(X,Y)=\sum_{(i,j)\in\nn^*\times \nn^*}a_{ij}X^{i}Y^j$. Alors cette série est inversible
et $S^{-1}(X,Y)=\sum_{k\in\nn}T^k$ est la série formelle associée à la fonction holomorphe $\alpha^{-1}$ sur $U^2$. Si on pose $T^k(X,Y)=
\sum_{(i,j)\in\nn^*\times \nn^*}a_{ij}^{(k)}X^{i}Y^j$, alors $ \{(i,j);\;a_{ij}^{(k)}\not= 0\}\subset \C_{+}$ car $\C_{+}$ est stable par addition. On en déduit d'une part que $\spec({\alpha^{-1}})\subset \C_{+}$ et d'autre part, par le théorème d'Abel, que $\spec({u^{-1}})\subset \C_{+}$. 
\end{preuve}
\begin{preuve}{du théorème \ref{clef}}
La fonction $f$ admet sur $ \tore^2$ un minimum strictement positif, puisqu'elle est semi-continue inférieurement et on peut quitte à normaliser supposer que $f>1$ ce qui permet alors d'affirmer que $\ln f=\ln_{+}f$ est semi-continue inférieurement  et que $\ln f\in L^\infty(\tore^2)\subset L^1(\tore^2)$. Par le lemme \ref{sing}, il existe une mesure positive $\sigma$ singulière sur $\tore^2$ et une fonction holomorphe $g$ sur $U^2$, telles que $P[\ln f-\sigma]=\re(g)$ et $\spec({g})\subset\C_{+}$. Puisque $\exp(g)\in H(U^2)$ et que $U^2$ est simplement connexe, on peut trouver une fonction $u\in H(U^2)$ telle que $u^2=\exp(g)$. 
 En fait $u$ appartient à $H^\infty(U^2)$. En effet $\re(g)<P[\ln f]<||\ln f||_{\infty}$ et $|u^2|=|u|^2=\exp\big(\re(g)\big)$. Notons alors $u^*$ sa limite radiale sur $\tore^2$. L'égalité $\re(g)=P[\ln f-d\sigma]$ entra\^{i}ne $\ln f=\re(g)^*$. D'où on déduit que $f=|u^*|^2=\alpha\bar{\alpha}$ avec $\alpha=u^*$.
De l'égalité  $u=\exp(\frac{1}{2}g)$, on déduit, avec le corollaire \ref{raide} et le lemme \ref{unsur}, que $\spec({\alpha^{\pm 1}})\subset\C_{+}$.
\end{preuve}
\begin{preuve}{du corollaire \ref{clebar}}
Soit $f$ vérifiant les hypothèses du théorème \ref{clef}. On peut toujours supposer, quitte à restreindre $\C=\C^+\bigcup-\C^+$ que $\C^+$ est engendré sur $\nn$ par une $\zz$-base $\{e_{1},e_{2}\}$. Soit alors $\{\varepsilon_{1},\varepsilon_{2}\}$ une $\zz$-base de $\zz^2$ incluse dans $Y_{2}^+$. On pose $\tilde{\ma C}^+$ le demi-cône engendré sur $\nn$ par $\{\varepsilon_{1},\varepsilon_{2}\}$. Il existe un unique
automorphisme de $\zz^2$, noté $s$, tel que $s(e_{i})=\varepsilon_{i},i=1,2$. \'{E}crivons maintenant $e_{1}=(a,b), e_{2}=(c,d),\varepsilon_{1}(\alpha,\beta),\varepsilon_{2}=(\gamma,\delta )$. Si $M=(u,v)$ est un point de $\ma{C}^+$ , alors $s\big((u,v)\big)=(u',v')$ avec $\left(\begin{matrix}u'\\v'\end{matrix}\right)=\left(\begin{matrix}\alpha&\gamma\\-\beta&\delta \end{matrix}\right)\left(\begin{matrix}d&-c\\-b&a\end{matrix}\right)\left(\begin{matrix}u\\v\end{matrix}\right)$, $\alpha\delta -\beta\gamma=ad-bc=1$.
On définit un morphisme $\tilde{s}$ de $\tore^2$ sur  lui-m\^{e}me en posant $\tilde{s}(\chi_{1}^{u} \chi_{2}^v)=\chi_{1}^{u'}\chi_{2}^{v'}$ étendu par linéarité. Alors la fonction $f\circ \tilde{s}$ vérifie les hypothèses du théorème \ref{clef} grâce auquel on peut écrire $f\circ \tilde{s}=\alpha\bar{\alpha}$, avec $\alpha\in H^{++}$. On en déduit que $f=(\alpha\circ {\tilde{s}}^{-1})\overline{(\alpha\circ {\tilde{s}}^{-1})}$ où $\alpha\circ {\tilde{s}}^{-1}\in H^{\C^+}$.
\end{preuve}

 \section{ Une application du théorème \ref{clef}}\label{polo}
 \subsection{ Factorisations du symbole associées au polygone}\label{factor}
 \subsubsection{Notations et définitions}\label{notar}
 Dans cette section, on définira l'opérateur de Toeplitz tronqué sur un polygone convexe, associé à un symbole régulier $f$ et on déduira du 
théorème \ref{clef} et de son corollaire (corollaire \ref{clebar}) un théorème  de \og factorisation minimale\fg \;du symbole (proposition \ref{terminator}), factorisation liée à la géométrie du polygone.

 Soit $\widetilde{\Lambda }$ un polygone de $\rr^2$ de sommets $\{A_{i}\}_{i=1,\ldots,m}$. Sans restreindre la généralité du problème on supposera que les coordonnées des sommets $A_{i}$ sont des \emph{entiers naturels}.  Au besoin $A_{m+1}$ désignera $A_{1}$ : en fait on a intér\^{e}t à considérer les indices des sommets dans $\zz/m\zz$.  La droite portant le côté $c_{i}=[A_{i},A_{i+1}]$ est notée $d_{i}$. On notera $\tilde S_i^+$ les demi-espaces définis par $d_{i}$ de sorte que $\widetilde\Lambda =\bigcap_{i=1}^n\tilde S_i^+$ et $ \tilde S_i^-=\rr^2\setminus \tilde S_i^+$. Le demi-espace $\tilde{S_{i}}^+,$ est appelé le \emph{demi-espace positif déterminé par la droite} $d_{i}$, le demi-espace $\tilde{S_{i}}^-$ étant le\emph{ demi-espace négatif déterminé par} $d_{i}$.  Soit $\tau_{u_{i}}$  la translation de vecteur $u_{i}=\overrightarrow{A_{i}O}$. Notons $\tau_{u_{i}}(d_{i})=d_{i,0}$ et $\tau_{u_{i}}(\tilde{S}_{i}^+)=\tilde{S}_{i,0}^+$.
Notre étude concerne les traces des espaces précédents sur $\zz^2$. Nous aurons ainsi à considérer 
$\Lambda =\widetilde{\Lambda }\cap \zz^2
 $, $S_k^+=\tilde S_k^+\cap \zz^2$,  $S_k^-=\tilde S_k^-\cap \zz^2$, $S_{k,0}^+=\tilde S_{k,0}^+\cap \zz^2$. On a donc 
$\Lambda=\bigcap_{k=1}^n S_k^+$.  
On note $\C_{A_{i}}$ le cône de sommet $A_{i}$ déterminé par les deux droites $d_{i-1}$ et $d_{i}$, $\C_{A_{i}}^+$ le demi-cône de 
$\C_{A_{i}}$ qui contient $\Lam$, $\C_{A_{i,0}}$ (cône déterminé par les deux droites $d_{i-1,0}$ et $d_{i,0}$) et $\C_{A_{i,0}}^+$ les 
translatés de vecteur $u_{i}$ de  respectivement   $\C_{A_{i}}$ et $\C_{A_{i}}^+$.

\begin{figure}[!]\centering\includegraphics[width=.3\textwidth]{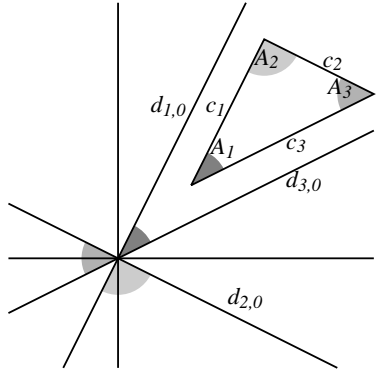}
 \caption{Le polygone}\label{trian}
 \end{figure}

La figure \ref{trian} illustre la description précédente dans le cas d'un triangle, les demi-cônes positifs étant grisés dans le triangle $A_{1}A_{2}A_{3}$.

 Pour tout $i\in\{1,\ldots,n\}$ on notera $p_{i}$ la projection orthogonale de $L^2(\tore^2)$ sur $\ma{P}(S_{i}^-)$.

Nous allons voir que sous les hypothèses décrites dans le théorème \ref{clef} sur $f$ et sous des hypothèses minimales sur $\Lambda$, l'opérateur $T_{\Lam}(f)$ est inversible et établir une expression de son inverse. Une des clefs de l'inversion des opérateurs de Toeplitz tronqués sur un polygone est la d\'ecomposition du symbole $f$ relativement \`a chaque face du polygone. Cette d\'ecomposition est l'objet du théorème suivant.
\begin{proposition}\label{zer}
On suppose  que la fonction $f$ vérifie les hypothèses du théorème \ref{clef}. Pour tout sommet $A$ du polygone convexe $\Lam$, il existe un cône  $C_{A}\subset \C_{A,0}$  pour lequel la décomposition $f=\alpha\bar{\alpha}$ du corollaire \ref{clebar} vérifie la propriété de compatibilité suivante :
pour tout $i\in\zz/m\zz$ le symbole $f$ admet une décomposition $f=G_{i}\ovl{G}_{i}$, avec
\begin{equation*}
G_{i}^{\pm 1}\in L^\infty(\tore^2),\;\; 
 \spec({G_{i}})\subset S_{i,0}^+\;\text{et}\;G_{i}\in\{\alpha,\bar{\alpha}\}.
 \end{equation*}
\end{proposition}
\begin{preuve}{de la proposition \ref{zer}}
Soit $A$ un sommet de $\Lam$ et avec les notations habituelles $\C_{A}=\C_{A}^+\cup\C_{A}^-$ le cône de sommet $A$ déterminé par les deux côtés du polygone issus de $A$ de sorte que $\C_{A}^+\cap\Lam=\Lam$.  Avec les notations du début de cette section, on considère le faisceau de droites en $O$ constitué par la famille $\mathcal{F}=\{d_{i,0}\}_{i=1,\ldots,m}$.
Pour un certain entier $p\in]0,n]$, un nombre $n-p\ge 2$ de demi-droites des droites de $\ma F$  sont incluses dans l'intérieur de $\C_{A,0}^+$,
constituant ainsi  une partition de $\C_{A,0}^+$ en $n-p-1$ demi-cônes 
 définis par ces demi-droites. Choisissons arbitrairement un de ces demi-cônes : il correspond à un cône $C_A$ de sommet $O$, inclus dans $\C_{A,0}$
défini par deux droites de la famille $\mathcal{F}$. Il n'existe pas de droite de $\mathcal{F}$
incluse dans le cône $C_{A}$.
Le  corollaire \ref{clebar} donne une factorisation $f=\alpha\bar\alpha$ avec $\spec(\alpha)\subset C_{A}^+\subset\C_{A,0}^+$. On peut alors écrire $f=G_{i}\ovl{G_{i}}$ avec 
$G_{i}= \alpha\;\text{si}\;C_{A}\subset {S_{i,0}}^+$ et $G_{i}=
\bar{\alpha}\;\text{si}\;C_{A}\subset {S_{i,0}}^-.
$
\end{preuve}
\begin{definition}\label{cancan}
Avec les notations de la proposition \ref{zer}, le $m$-uplet $(G_{1},\ldots,G_{m})$ est une factorisation  de $f$ associée à $(\Lam,A,{C}_{A})$.
\end{definition}

Par exemple, si $\Lam$ désigne le triangle représenté par la figure \ref{trian}, et si on choisit pour chaque sommet $A_{i}$ le cône  $\C_{{i}}=\mathcal{C}_{A_{i},0}$,\ et si on note $f=\alpha_{i}\bar{\alpha}_{i}$ la décomposition de $f$ par rapport au cône $\mathcal{C}_{i}$ prévue par le théorème \ref{clef} et son corollaire, avec $\spec(\alpha_{i})\subset \C_{i}^+$,  le triplet $(\alpha_{1},\bar{\alpha}_{1},\alpha_{1})$ (respectivement $(\alpha_{2},{\alpha}_{2},\bar{\alpha}_{2})$,$(\bar{\alpha}_{3},{\alpha}_{3},\alpha_{3})$) est une
 factorisation  de $f$ associée à $(\Lam,A_{1},\mathcal{C}_{1})$ (respectivement à
  $(\Lam,A_{2},\mathcal{C}_{2})$  et  $(\Lam,A_{3},\mathcal{C}_{3})$).  
On peut faire deux remarques. 
\begin{enumerate}
\item Une factorisation associée à un triplet $(\Lam,A,C_{A})$ est liée à la numérotation des c\^{o}tés de $\Lam$ et à la factorisation du symbole dans le cône $C_{A}\subset \C_{A,0}$. \`{A} ceci près il y a unicité.
\item Si on considère pour un sommet $A$ du polygone $\Lam$ à $m$ sommets, une factorisation $\mathbb{F}$  associée à $(\Lam,A,\ma{C}_{A})$, il existe pour un entier positif $p\ge 2$ une suite $\{m_{0}=0<1= m_{1}< m_{2}<\ldots< m_{p}=m\}$ définissant la partition de $\{1,\ldots,m\}=\bigcup_{i=1}^pI_{i}$ où $I_{i}=\{m_{i-1}+1,\ldots, m_{i}\}$,  vérifiant la propriété suivante :
\begin{equation}\label{smurt}
\forall i\in\{1,\ldots,p\}\;\forall k,l\in I_{i}\; \;\;G_{k}=G_{l}\equiv g_{i}, g_{i}\not= g_{i+1}\;\text{si on pose } \;g_{p+1}= g_{1}
\end{equation}

Pour tout $i\in\{1,\ldots,p\}$, notons \begin{equation}\label{seres}
\So_{i}^-=\bigcup_{k\in I_{i}}S_{k}^-\;\;\So_{i}^+=\zz^2\setminus\So_{i}^-\;\;\text{et}\;\;\So_{i,0}^+=\bigcap_{I_{i}}S_{i,0}^+.
\end{equation}
Et alors on a :
\begin{enumerate}
\item  $\Lambda=\bigcap_{i=1}^p{\So_{i}^+}$
\item Pour tout $i\in \{1,\ldots,p\}$, on a  $\spec ({{g}_{i}})\subset \So_{i,0}^+$ et ${{g}_{i}}\in\{\alpha, \bar{\alpha}\}$.
\item $ \Po(\Lam)=\bigcap_{i=1}^p{\Po(\So_{i}^+)}.$
\end{enumerate}
\end{enumerate}
Par exemple, si on renumérote les sommets de la figure \ref{trian} selon $B_{1}=A_{3},B_{2}=A_{1},B_{3}=A_{2}$, et si on considère la factorisation du symbole $f=\alpha_{3}\overline{\alpha_{3}}$ relativement au cône $\ma C_{3}=\ma C_{B_{1},0}$, la factorisation  associée à  $(\Lam, B_{1},\ma C_{3})$ est $(\overline{\alpha_{3}},\alpha_{3},\alpha_{3})$ et on a $I_{1}=\{1\},I_{2}=\{2,3\}, g_{1}=\overline{\alpha_{3}},g_{2}=\alpha_{3}$.

\begin{definition}\label{derive}
Avec les notations de la proposition \ref{zer}, le $p$-uplet $(g_{1},\ldots,g_{p})$, $p\ge 2$ est la \emph{factorisation de $f$ définie par la somme hilbertienne $F=\bigoplus_{i=1\ldots,p}^h\Po(\So_{i}^-)$, issue de la factorisation $\mathbb{F}$ associée à $(\Lam,A,\ma{C})$}. Si $p=2$, on dit que $(g_{1},g_{2})$ est une factorisation minimale.
\end{definition}
L'intér\^{e}t de la factorisation $(g_{1},\ldots,g_{p})$ est d'utiliser moins de demi-espaces que $\mathbb F$ pour reconstituer $\Lam$ par intersection . Cet intér\^{e}t ainsi que
le sens de la somme hilbertienne $F$ introduite dans la définition \ref{derive} appara\^{i}tront dans le théorème d'inversion en sous-section \ref{invoptotro}.

\begin{figure}[!]\centering\includegraphics[width=.3\textwidth]{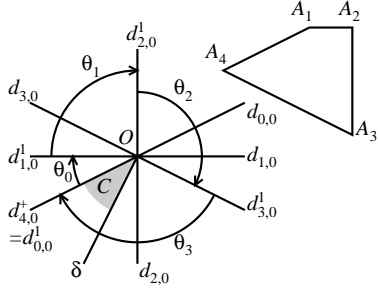}
 \caption{Factorisation minimale, n=$4$.}\label{facmin}
 \end{figure}

La proposition suivante assure l'existence de factorisations minimales.
\begin{proposition}\label{terminator}
 On suppose  que la fonction $f$ vérifie les hypothèses du théorème \ref{clef}. Pour tout sommet $A$ de $\Lam$, il existe un cône ${C}\subset
\ma{C}_{A,0}$ pour lequel existe une factorisation  minimale issue de la factorisation canonique associée à $(\Lam,A,{C})$.
\end{proposition}
\begin{preuve}{de la proposition \ref{terminator}}
Une construction possible de la factorisation minimale est la suivante.
 On numérote les sommets à partir de $A$ dans le sens des aiguilles d'une montre (sens dit ici positif). On a donc $A=A_{1}$. En reprenant les notations du paragraphe \ref{notar},
 soit $\ma{F}_{0}=\{d_{1,0},\ldots,d_{m,0}=d_{0,0}\}$ le faisceau de droites passant par $O$ déterminé par $\Lam$. On définit  $d_{m,0}^+\equiv d_{0,0}^1$ comme étant la translatée de vecteur $\overrightarrow{A_{1}O}$ de la demi-droite de $d_{m}^+$ telle que $\ma C_{A_{1}}^+= \ma C_{A_{1}}^+[d_{m}^+,d_{1}^+]$ (voir notation $7$ de la section \ref{nonotes}). Pour tout $j\in\{0,\ldots,m\}$, on définit par récurrence la demi-droite $d^1_{j+1,0}$  comme la première demi-droite de la première droite $d_{j+1,0}$ du faisceau $\ma{F}_{0}$ rencontrée en tournant dans le sens positif à partir de $d^1_{j,0}$ et on pose   $\theta_{j}\in[0,2\pi]$  la mesure en radians de l'angle orienté positivement formé par $d_{j,0}^1$ et  $d_{j+1,0}^1$ et $\theta_{m}=\theta_{0}$. Soit $r\in\nn$ défini par 
 \[\theta_{0}+\theta_{1}+\ldots+\theta_{r}<\pi\le \theta_{0}+\theta_{1}+\ldots+\theta_{r}+\theta_{r+1}.\]
On pose $S_{n}^-=S_{0}^-$,  et   $\So_{1}^-=\bigcup_{i\in\{0,1,2,\ldots,r+1\}}S_{i}^-, \;\So_{2}^-=\bigcup_{i\notin\{0,1,2,\ldots,r+1\}}S_{i}^-$.
 On définit ${C} \subset {C}_{A,0}$ par  ${C}={C}[d_{m,0},\delta ]$ où $\delta $ est une droite passant par $O$ incluse dans $C_{A,0}$ de sorte que le cône ${C}$ ne contienne aucune droite de $\ma{F}_{0}$. Le théorème \ref{clef} assure une décomposition de $f$ sous la forme
 $f=\alpha\bar{\alpha}$, avec $\spec (\alpha)\subset {C}^+$. Tout a été fait pour que $(\alpha,\bar{\alpha})$ soit la factorisation minimale annoncée définie par la somme hilbertienne $F=\ma{P}(\So_{1}^-)\bigoplus \ma{P}(\So_{2}^-)$.
 \end{preuve}
 Cette construction est illustrée par la figure \ref{facmin} si l'on prend $A=A_{1}, \theta_{0}+\theta_{1}<\pi< \theta_{0}+\theta_{1}+\theta_{2}+\theta_{3}, \So ^-_{1}=S_{0}^-\cup S_{1}^-\cup S_{2}^- $et $\So_{2}^-=S_{3}^-$.

\subsection{Inversion de l'opérateur de Toeplitz tronqué sur un polygone convexe.}\label{invoptotro}
\emph{Dans ce qui suit, nous supposerons que $f$ est une fonction positive, essentiellement bornée et semi-continue inférieurement sur le tore $\tore^2$.}
Soit $\Lam$ un polygone convexe, $A$ un sommet et $\mathbb F$ une factorisation associée à un triplet 
$(\Lam,A,C_{A})$ (voir définition \ref{cancan}). On se donne conformément à la définition \ref{derive} une factorisation issue de $\mathbb F$ notée $(g_{1},\ldots,g_{p})$ définie par l'espace
 \begin{equation}\label{hilb}F=\bigoplus_{i=1\ldots,p}^h\Po(\So_{i}^-),\end{equation} 
appelé \emph{l'espace fondamental de l'inversion}
(la lettre $h$ pour rappeler qu'il s'agit d'une somme hilbertienne).

Cette factorisation  définit une famille d'\emph{op\'erateurs de Hankel} $H(i,j)$   comme suit.
\begin{definition}
\label{def03}
Posons pour tout $i\in\{1,\ldots,p\},\;\Phi_{i}=\dfrac{\bar{{g}}_{i}}{\bar{{g}}_{1}}$ et d\'esignons par $\Pi_{i}$ la projection orthogonale de $L^2(\tore^2)$ sur $\Po(\So_{i}^-)$. Alors pour tout couple $(i,j) \in \left\{ 1,\;\ldots ,\;p\right\}^2$,
 on d\'efinit   un \emph{op\'erateur d'\'echange } $H(i,j)$ par :
\[
\begin{array}{rcl}
H(i,j):\;\Po(\So_j^-)& \longrightarrow &\Po(\So_i^-)\\
H(i,j)\theta_j &=& \Pi_i\left( \frac{\Phi_j}{\Phi _i}\theta_j\right)
\end{array}
\]

\end{definition}
Ces op\'erateurs de Hankel d\'efinissent eux-m\^{e}me un endomorphisme de $F$, not\'e $H$, de la mani\`ere suivante. Soit $\theta=(\theta_{1},\ldots,\theta_{m})$ un vecteur de $F$. Alors 
\begin{equation}\label{lebo}
H\theta=\big((H\theta)_{i}\big)_{i=1,\ldots,p} \;\text{o\`u}\;(H\theta)_{i}=\sum_{j=1}^pH(i,j)\theta_{j}
\end{equation}
\begin{proposition}\label{herm}
L'op\'erateur $H$ est un op\'erateur hermitien positif sur $F$.
\end{proposition}
\begin{preuve}{de la proposition \ref{herm}}
Soit $\theta$ et $\mu$ deux vecteurs de $F$, $\theta=(\theta_{1},\ldots,\theta_{p}),\mu=(\mu_{1},\ldots,\mu_{p})$. Par définition de la somme hilbertienne $F$, on a   $\langle \theta,\mu\rangle=\sum_{i=1}^p\la \theta_{i},\mu_{i}\ra_{\Po(\So_{i}^-)}$. Le produit scalaire de $\Po(\So_{i}^-)$  est en fait la restriction sur cet espace du produit scalaire sur $L^2(\tore^2)$. On le notera plus simplement $\la,\ra$. Alors en tenant compte du fait que $\Pi_{i}$ est autoadjoint et que $\bar{\Phi}_{i}=\dfrac{1}{\Phi_{i}}$ (à partir de la définition des $g_{i}$ dans la proposition \ref{zer}), on a :
\begin{align*}
&\langle H\theta,\mu\rangle_{F}=\sum_{i=1}^p\left\langle \sum_{j=1}^pH(i,j)\theta_{j},\mu_{i}\right\rangle=\sum_{i=1}^p\left\langle \sum_{j=1}^p\Pi_i\left( \frac{\Phi_j}{\Phi _i}\theta_j\right),\mu_{i}\right\rangle\\
&= \sum_{i=1}^p\left\langle \sum_{j=1}^p\theta_{j},\Pi_{j}\left(\frac{\Phi_{i}}{\Phi_{j}}\mu_{i}\right)\right\rangle=\sum_{j=1}^p\left\langle \sum_{i=1}^p\theta_{j},H(j,i)\mu_{i}\right\rangle=\langle \theta,H\mu\rangle_{F},
\end{align*}
ce qui prouve le caractère hermitien de $H$. Montrons sa positivité. On a pour tout $\theta\in F$ :

$\la H\theta,\theta\ra_{F}=\sum_{i=1}^p\left\langle \sum_{j=1}^p \Pi_i\left( \frac{\Phi_j}{\Phi _i}\theta_j\right),\theta_{i}\right\rangle=
\sum_{i=1}^p\left\la\Pi_i\Big( \frac{1}{\Phi _i} \sum_{j=1}^p \Phi_j\theta_j\Big),\theta_{i}\right\ra=$

\hspace{1,65cm}$=\sum_{i=1}^p\left\la \sum_{j=1}^p \Phi_j\theta_j,\Phi_{i}\theta_{i}\right\ra=\left\la \sum_{j=1}^p \Phi_j\theta_j,\sum_{i=1}^p\Phi_{i}\theta_{i}\right\ra$ et on conclut.
 \end{preuve}
 \begin{coro}\label{neuneu}
 La restriction $H_{\left|\im H\right.}$ de $H$ à $\im H$ est un automorphisme de $\im H$. 
 On a de plus $H_{|\im H}=I-\mathcal{M}_{}$ où $\mathcal{M}_{}$
 est défini  pour tout $\theta\in \im H$ par $(\mathcal{M}\theta)_{i}=\sum_{j=1}^pM(i,j)\theta_{j}$ avec $$M(i,j)=\begin{cases} -H(i,j)\;\text{si}\;j\not=i\\0\;\text{si}\;j= i\end{cases}$$
 \end{coro}
  \begin{preuve}{du corollaire \ref{neuneu}}		
 La restriction de $H$ à $\im H$ est un opérateur de $\im H$, inversible puisque $H$ étant hermitien $\ker H\cap\im H=\{0\}$. 
  \end{preuve}

\begin{definition}\label{def04} \hspace{10cm}
\begin{enumerate}
\item On appelle \emph{ champ de vecteurs fondamental associé à} la factorisation $(g_{1},\ldots, g_{p})$ l'application $\gamma(q)$ de $\Po(\Lam)$ dans $F$ définie pour tout polynôme $q\in\Po(\Lam)$ par   $\gamma(q)=(\gamma_{i}(q))_{i=1,\ldots,p}$ avec $\gamma_{i}(q)=\Pi_{i}(\bar{\Phi}_{i}\frac{q}{\bar{{g}}_{1}})$.
\item L'équation d'inconnue $\theta\in F$, $H\theta=\gamma$ est \emph{l'équation de Hankel} associée à la factorisation $(g_{1},\ldots,g_{p})$.
 \end{enumerate}
\end{definition}
\begin{theorem}\label{invon}
Soit $f$ une fonction vérifiant les hypothèses du théorème \ref{clef}, $\Lam$ l'intersection avec $\zz^2$ d'un polygone convexe de $\rr^2$ à $n$ sommets et une factorisation $(g_{1},\ldots,g_{p})$ obtenue conformément à la définition \ref{derive}.
  Alors 
\begin{enumerate}
\item L'équation de Hankel associée à $(g_{1},\ldots,g_{p})$  admet une unique solution dans $\im H$.
\item $T_{\Lam}(f)$ est inversible et avec les notations de la définition \ref{def04}, on a 
\begin{equation}\label{plouk}
T_{\Lam}(f)^{-1}(q)=\frac{q}{f}-\frac{1}{g_{1}}\sum_{i=1}^p\Phi_{i}\theta_{i}
\end{equation}
où $\theta$ est l'unique solution de \emph{l'équation de Hankel} $H\theta=\gamma(q)$ dans $\im H$.
\end{enumerate}
\end{theorem}
\subsection{Démonstration du théorème \ref{invon}}\label{demimpocon}

La clef de la formule d'inversion (\ref{plouk})  est le lemme \ref{Lemme010} : il exprime l'inverse de l'op\'erateur de Toeplitz tronqu\'e par l'interm\'ediaire d'une projection sur un sous-espace vectoriel ferm\'e de $L^2(\tore^2)$.
Posons \begin{equation}\label{ka}
K=\bigcap_{i=1}^p\Phi_{i}{\Po(S_{i}^+)}.
\end{equation} et $P_{K}$ la projection orthogonale de $L^2(\tore^2)$ sur $K$.                                                                                                        

\begin{lem}\label{coro007}
on pose $K_0=f \Po(\Lambda)$. Alors $K_{0}=\bar{g}_{1}K.$
\end{lem}
\begin{preuve}{du lemme \ref{coro007}}
Si $x\in \bar{g}_{1}K$, alors pour tout $i\in\{1,\ldots,p\}$, on a $x\in \bar{g}_{i}{\Po(S_{i}^+)}=f\left(g_{i}^{-1}{\Po(S_{i}^+)}\right)$.
Or  $\spec(g_{i}^{-1})\subset S_{i,0}^+$, ce qui implique pour tout $i$ l'inclusion $g_{i}^{-1}{\Po(S_{i}^+)}\subset 
{\Po(S_{i}^+})$. Ainsi $x\in \bigcap_{i=1}^pf\Po(S_{i}^+)=f\Po(\Lambda)$. 
 Réciproquement si $x\in f\Po(\Lam)=\cap_{i=1}^pf\Po(S_{i}^+)$, il existe 
pour tout $i\in\{1,\ldots,p\}$ un élément $y_{i}\in \Po(S_{i}^+)$ tel que $x=\bar{g}_{i}(g_{i}y_{i})$ et $x/\bar{g}_{1}=\Phi_{i}(g_{i}y_{i})\in\Phi_{i}\Po(S_{i}^+)$ d'où $x\in\bar{g}_{1}K$.
\end{preuve}
\begin{lem}\label{Lemme010} 
L'op\'erateur de Toeplitz $T_{\Lambda}(f)$ est inversible. De fa\c con plus pr\'ecise, si
$P_{K_0}$ est  la projection orthogonale de
$L_{1/f}^2(\tore^2)$ dans $K_0$.
Alors
\[
T_{\Lambda}^{-1}(f)(q)=\frac{1}{f}P_{K_0}(q) 
\]
\end{lem}
\begin{preuve}{du lemme \ref{Lemme010}}
Soit $q\in \Po(\Lambda)$. Il existe $p\in \Po(\Lambda)$ tel que 
$
P_{K_0}(q)=fp. 
$
Or $q=fp+(q-fp)$ et $\Pi_\Lambda(q-fp)=0$. En effet, pour tout $\delta \in \Po(\Lambda),$ on a par d\'efinition de $P_{K_{0}}$ l'égalité  $\langle q-P_{K_{0}}(q),f\delta \rangle_{1/f}=0$  qui s'\'ecrit encore 
$\langle q-fp,\delta \rangle=0$.
D'o\`u $q=\Pi_{\Lambda }(fp)=T_\Lambda(f)(p)$.
On en d\'eduit que $T_\Lambda(f)$ est une bijection de $\Po(\Lambda)$ sur
lui-m\^eme et 
$
T_\Lambda^{-1}(f)=\frac{1}{f}P_{K_0}. 
$
\end{preuve}
\begin{proposition}\label{proj}
L'opérateur $T_{\Lambda}(f)$ est inversible et pour tout polynôme $q\in \Po(\Lambda)$, on a :
\begin{equation}
T_{\Lambda}(f)^{-1}(q)=\frac{1}{g_{1}}P_{K}(\frac{q}{\bar{g}_1})
\end{equation}
\end{proposition}
\begin{preuve}{de la proposition \ref{proj}}
Pour tout polynôme $q\in \Po(\Lambda)$, on v\'erifie que  $P_{K}(\frac{q}{\bar{g}_{1}})=\frac{1}{\bar{g}_{1}}P_{K_{0}}(q)$.
En effet, si $\delta \in K$ et $q\in \Po(\Lambda)$, on a l'\'egalit\'e $0=\langle P_{K}(\frac{q}{\bar{g}_{1}})-\frac{q}{\bar{g}_{1}},\delta \rangle=\langle \bar{g}_{1}P_{K}(\frac{q}{\bar{g}_{1}})-q,\bar{g}_{1}\delta \rangle_{1/f}.$ Les lemmes \ref{coro007} et \ref{Lemme010} permettent de conclure.
\end{preuve}
La proposition suivante donne acc\`es au calcul explicite du projecteur $P_{K}$ et par suite de l'inverse de l'op\'erateur de Toeplitz tronqu\'e.
\begin{proposition} \label{ptibig}
Soit $\psi\in L^2(\tore^2)$. On pose $\Gamma(\psi)\in F$ le vecteur $(\Gamma_{i})_{i=1\ldots,p}$ o\`u $\Gamma_{i}=\Pi_{i}(\bar{\Phi}_{i}\psi)$.
L'\'equation $H\theta=\Gamma(\psi)$ d'inconnue $\theta$ admet une unique solution dans $\im\;H$. Si $(\theta_{1},\ldots,\theta_{p})$ est la d\'ecomposition de cette solution dans $F$, alors\[P_{K}(\psi)=\psi-\sum_{i=1}^p\Phi_{i}\theta_{i}.\]
\end{proposition}
La d\'emonstration de la proposition \ref{ptibig} est la cons\'equence imm\'ediate des trois lemmes suivants.
\begin{lem}\label{ach}
Le noyau de $H$ est l'ensemble $
\{(\alpha_{1},\ldots,\alpha_{p})\in F; \sum_{i=1}^p\alpha_{i}\Phi_{i}=0\}$.
\end{lem}
\begin{lem}\label{achach}
Pour tout vecteur  $\psi$ de $L^2(\tore^2)$, il existe un unique vecteur $\theta$ de $\im\;H$ tel que $P_{K}(\psi)=\psi-\sum_{i=1}^p\Phi_{i}\theta_{i}.$
 \end{lem}
\begin{lem}\label{achachach}
 Avec les notations de la proposition \ref{ptibig}, il y a \'equivalence entre 
\begin{itemize}
\item[i)] Le vecteur $\theta$ de $F$ est solution de $H\theta=\Gamma(\psi)$.
\item[ii)] $P_{K}(\psi)=\psi-\sum_{i=1}^p\Phi_{i}\theta_{i}$.
\end{itemize}
\end{lem}
Prouvons les lemmes \ref{ach}, \ref{achach}, \ref{achachach}.
\begin{preuve}{du lemme \ref{ach}}
Si l'on suppose $\alpha\in\ker H$, alors $\forall i\in\{1,\dots,p\}\;\sum_{ j=1}^p\Pi_{i}(\frac{\Phi_{j}}{\Phi_{i}}\alpha_{j})=0.$
Comme $\alpha_{i}\in \Po(S_{i}^-)$ et $\ker \Pi_{i}={\Po(S_{i}^+)}$, cette égalité traduit d'une part que $\sum_{j=1}^p\Phi_{j}\alpha_{j}\in  \bigcap_{i=1}^p \Phi_{i}{\Po(S_{i}^+)}=K$ et $\sum_{j=1}^p\Phi_{j}\alpha_{j}\in \sum_{j=1}^p\Phi_{j}\Po(S_{j}^-)=K^\perp$. D'où
$\sum_{j=1}^p\Phi_{j}\alpha_{j}=0$.
La r\'eciproque est imm\'ediate.
\end{preuve}
\begin{preuve}{du lemme \ref{achach}}
Pour l'unicit\'e, remarquons que si $\theta\in\im\;H$ est tel que $\sum_{i=1}^p\theta_{i}\Phi_{i}=0$ alors $\theta\in \im\;H\cap\ker  H=\{0\}$
puisque $H$ est hermitien (proposition \ref{herm}). L'existence r\'esulte du fait que 
$K^\perp=\sum_{i=1}^p \Phi_{i}\Po(S_{i}^-)$, ce qui, pour tout $\psi\in L^2(\tore^2)$,  permet d'\'ecrire  la projection orthogonale $P_{K^\perp}(\psi)$ sous la forme $\sum \theta_{i}\Phi_{i}$ avec $\theta\in F$ et on conclut, pour montrer qu'on peut prendre $\theta$ dans $\im H$, en  d\'ecomposant  $F$ en la somme directe $F=\ker H\bigoplus \im\;H$ et en utilisant le lemme \ref{ach}.
\end{preuve}
\begin{preuve}{du lemme \ref{achachach}}
Supposons $i)$. Alors $\forall i\in\{1,\dots,p\}\;\Pi_{i}\big(\bar{\Phi}_{i}(\sum_{k=1}^p\theta_{k}\Phi_{k}-\psi)\big)=0$. En n'oubliant pas que 
$\Phi_{i}\bar{\Phi}_{i}=1$, on en d\'eduit que $\sum_{k=1}^p\theta_{k}\Phi_{k}-\psi\in K$ et du fait que $\sum_{k=1}^p\theta_{k}\Phi_{k}\in K^{\perp}$, on a  $\sum_{k=1}^p\theta_{k}\Phi_{k}=
P_{K^\perp}(\psi)$ et par cons\'equent $P_{K}(\psi)=\psi-\sum_{i=1}^p\theta_{i}\Phi_{i}.$

\noindent Supposons $ii)$.  Pour tout $i$ on a $\Pi_{i}\Big(-\bar{\Phi}_{i}\psi+\sum_{j=1}^p\frac{\Phi_{j}}{\Phi_{i}}\theta_{j}\Big)=\Pi_{i}\big(\bar{\Phi}_{i}P_{K}(\psi)\big)=0$, ce qui implique $\Gamma_{i}=\Pi_{i}(\bar{\Phi}_{i}\psi)=(H\theta)_{i}$. Et on conclut.
\end{preuve}
\begin{preuve}{du théorème \ref{invon}}
Le deuxième item est la conséquence des  propositions \ref{proj} et \ref{ptibig}. Le premier item découle directement des lemmes \ref{ach} et \ref{achach} qui montrent de plus que $\gamma(q)\in\im H$. Mais, selon le corollaire \ref{neuneu}, $\theta=H_{|\im H}^{-1}(\gamma(q)).$
\end{preuve}
Le théorème \ref{invon} donne une formule théorique d'inversion de $T_{\Lam}(f)$ qui passe par la résolution de l'équation d'inconnue $\theta(q)$, $H_{|\im H}(\theta)=\gamma(q)$.
Ceci est appliqué dans la section suivante lorsque $\Lam$ est un triangle.
\section{Le théorème de trace. Cas du triangle.}\label{lap}
 Si $\Lam_{\lambda}$ est le triangle décrite en sous-section \ref{trisym}, on peut écrire
conformément aux notations \ref{notar} $\Lam_{\lambda}=\bigcap_{i=1}^3S_{i}^+$. On pose $\ma{C}$ le cône formé des deux droites portant les côtés $c_{1}$ et $c_{3}$ du triangle, $\ma{C}^+=\ma{C}\cap\nn^2$. D'après le corollaire \ref{clebar}, $f$ se factorise suivant $f=\alpha\bar{\alpha}$ avec $\alpha\in H^{\ma{C}^+}$. Si on pose $\mathbb{S}_{1}^-=S_{1}^-\cup S_{3}^-,\mathbb{S}_{2}^-=S_{2}^-, \Pi_{1},\Pi_{2} $ les projections orthogonales de $L^2(\tore^2)$ sur respectivement $\mathcal{P}(\mathbb{S}_{1}^-)$ et $\mathcal{P}(\mathbb{S}_{2}^-)$, alors $(\alpha,\bar{\alpha})$ est la factorisation minimale issue de la factorisation $(\alpha,\bar\alpha,\alpha)$ associée à $(\Lam,O,\C)$,
 définie par  la somme hilbertienne $F=\ma{P}(\mathbb{S}_{1}^-)\bigoplus \ma{P}(\mathbb{S}_{2}^-)$ (voir les définitions \ref{cancan} et \ref{derive}).
Rappelons l'énoncé du théorème de trace( théorème \ref{otlitchno}).

\noindent \textbf{Théorème de trace}

\noindent On suppose que $f$ vérifie les hypothèses du théorème \ref{clef} et que la fonction $\frac{1}{f}$ est essentiellement bornée sur $\tore^2$.
Alors 
\begin{align*}
\trace\left(T_{\Lam_{\lambda}}(f)^{-1}\right)=|\Lam| \;||\frac{1}{f}||_{1}&+\left(\sum_{(u,v)\in\zz^2}|u|\;\widehat{\log \left(\frac{1}{f}\right)}(u,v)\overline{\widehat{\frac{1}{f}}(u,v)}\right)\mathfrak{S}_{1}(\lambda)\\&+\left(\sum_{(u,v)\in\zz^2}|v|\;\widehat{\log \left(\frac{1}{f}\right)}(u,v)\overline{\widehat{\frac{1}{f}}(u,v)}\right)\mathfrak{S}_{2}(\lambda)+o(\lambda).
\end{align*}

\subsection{Démonstration du théorème de trace}\label{trui}
\begin{lem}\label{inversons}
On considère la factorisation   $f=\alpha\overline{\alpha}$ décrite  en préambule du paragraphe \ref{lap}.
Posons
 \begin{align*}
 &\text{Pour tout}\; (k,l)\in\Lam_{\lambda},\;q=\chi_{1}^{k}\chi_{2}^l,\\
&\text{Pour tout}\; (k,l)\in\Lam_{\lambda},\;\zeta(q)=\Pi_{2}\big(\frac{q}{\alpha}\big)-\Pi_{2}\left(\frac{\overline{\alpha}}{\alpha}\Pi_{1}\big(\myfrac{q}{\overline{\alpha}}\big)\right),\\
  &I \text{ l'identité de}\; \ma{P}(\So_{2}^-),\\ 
  & \ma{N}_{}\;\text{ l'opérateur de}\; \ma{P}(\So_{2}^-)\;\text{ défini par l'égalité}\;
 \ma{N}(\theta_{2})=\Pi_{2}\left(\frac{\overline{\alpha}}{\alpha}\Pi_{1}\big(\myfrac{\alpha}{\overline{\alpha}}\theta_{2}\big)\right).
\end{align*}
 La résolution de l'équation de Hankel donne la formule d'inversion suivante pour tout  $(k,l)\in\Lam_{\lambda}$.
\begin{multline}\label{inversao}
 T_{\Lam_{\lambda}}(f)^{-1}(\chi_{1}^{k}\chi_{2}^l)=\\
\frac{\chi_{1}^{k}\chi_{2}^l}{f}-\frac{1}{\alpha_{}}\Pi_{1}\left(\frac{\chi_{1}^{k}\chi_{2}^l}{\bar{\alpha}}\right)-\frac{1}{\bar{\alpha}}(I-\ma{N})^{-1} \zeta(\chi_{1}^{k}\chi_{2}^l)+\frac{1}{\alpha}\Pi_{1}\left(\frac{\alpha_{}}{\bar{\alpha}}\big(I-\ma{N}\big)^{-1} \zeta(\chi_{1}^{k}\chi_{2}^l)\right).
\end{multline}
\end{lem}
\begin{preuve}{}
En se rappelant que $(\theta_{1},\theta_{2})$ désigne un élément de $F$, espace défini par l'équation (\ref{hilb}),
 l'équation de Hankel équivaut au système suivant  :
 \begin{align}[left=\empheqlbrace]
 &\theta_{1}+\Pi_{1}\big(\myfrac{\alpha}{\overline{\alpha}}\theta_{2}\big)=\Pi_{1}\big(\myfrac{q}{\overline{\alpha}}\big)\label{op}\\
&\Pi_{2}\big(\frac{\overline{\alpha}}{\alpha}\theta_{1}\big)+ \theta_{2}=\Pi_{2}\big(\frac{q}{\alpha}\big)\label{iop}
 \end{align}
 
 l'élimination de $\theta_{1}$ entre les équations (\ref{op}) et (\ref{iop})  donne directement l'équation
\begin{equation}\label{opi}
(I-\ma{N})(\theta_{2}(q))=\zeta(q)
\end{equation}
La démonstration s'achève directement avec le théorème \ref{invon} et l'égalité (\ref{plouk}).
\end{preuve}
 Rappelons que si $\phi$ est un opérateur de $L^2(\tore^2)$, sa norme de \emph{Hilbert-Schmidt associée au triangle }
 $\Lam$ se définit par $||\phi||_{\Lam}=(\sum_{(k,l)\in\Lam}||\phi(\chi_{1}^{k}\chi_{2}^l)||^2)^{1/2}$ (voir \cite{Z} ou \cite{DUSC}).
 \begin{coro}\label{tras}
 Avec les notations du lemme \ref{inversons}, posons $\ma W_{}(q)=-\Pi_{2}\big(\frac{\bar\alpha}{\alpha}\Pi_{1}(\frac{q}{\bar\alpha})\big),\ma V_{}(q)=\pi_{2}(\frac{q}{\alpha})$.
 La trace est donnée par l'égalité 
 \begin{equation}\label{trasa}
 \trace\left( T_{\Lam_{\lambda}}(f)^{-1}\right)=|\Lam_{\lambda|}\left|\left|\frac{1}{f}\right|\right|_{1}-\sum_{(k,l)\in\Lam_{\lambda}}\left|\left|\Pi_{1}(\frac{\chi_{1}^{k}\chi_{2}^l}{\bar\alpha})\right|\right|_{2}^2-\sum_{(k,l)\in\Lam_{\lambda}}\left|\left|\Pi_{2}(\frac{\chi_{1}^{k}\chi_{2}^l}{\alpha})\right|\right|_{2}^2 +R_{1}+R_{2}
 \end{equation}
 
 avec $|R_{1}|\le ||\ma W||_{\Lam}(||\ma W||_{\Lam}+2||\ma V||_{\Lam})$ et 
$ |R_{2}|\le C \sum_{(k,l)\in\Lam_{\lambda}}||\Pi_{1}(\frac{\alpha}{\bar\alpha}\zeta(\chi_{1}^{k}\chi_{2}^l))||_{2}^2$,
où $C$ est une constante.
 \end{coro}
 \begin{preuve}{du corollaire \ref{tras}}
 Par définition, $ \trace\left( T_{\Lam_{\lambda}}(f)^{-1}\right)=\sum_{\Lam_{\lambda}}\la T_{\Lam_{\lambda}}(f)^{-1}(\chi_{1}^{k}\chi_{2}^l),\chi_{1}^{k}\chi_{2}^l\ra$, ce qui donne immédiatement à partir de l'équation (\ref{inversao}) en posant 
 $\ma R(q)=-\frac{1}{\overline\alpha}(I-\ma N)^{-1}\zeta+\frac{1}{\alpha}\Pi_{1}(\frac{\alpha}{\overline\alpha}(I-\ma N)^{-1}\zeta)$ :
 
 \[ \trace\left( T_{\Lam_{\lambda}}(f)^{-1}\right)=|\Lam_{\lambda|}||\frac{1}{f}||_{1}-\sum_{(k,l)\in\Lam_{\lambda}}||\Pi_{1}(\frac{\chi_{1}^{k}\chi_{2}^l}{\overline\alpha})||_{2}^2+\sum_{(k,l)\in\Lam}\la\ma R , \chi_{1}^{k}\chi_{2}^l)\ra.\] Mais si on réécrit $\ma R (q)=-\frac{1}{\alpha}\Pi_{1}^\perp\big(\frac{\alpha}{\overline\alpha}(I-\ma N)^{-1}\zeta(q)\big)$  l'opérateur  $\Pi_{1}^\perp$désignant $I-\Pi_{1}$, on a
 
   $\sum_{(k,l)\in\Lam}\la\ma R , \chi_{1}^{k}\chi_{2}^l)\ra=
 -\sum_{(k,l)\in\Lam}\la (I-\ma N)^{-1}\zeta,\zeta\ra$.  Notant l'identité
 \begin{equation}\label{fuchs}
  (I-\ma N)^{-1}=I+\ma N(I-\ma N)^{-1},
  \end{equation}
   on a finalement en posant $R_{1}=||\ma W-\ma V||^2_{\Lam}-||\ma V||^2_{\Lam}$ et $R_{2}=\sum_{\Lam}\la (I-\ma N)^{-1}\zeta,\ma N\zeta\ra$ :
   
  \noindent$\sum_{(k,l)\in\Lam}\la\ma R , \chi_{1}^{k}\chi_{2}^l)\ra=-\sum_{(k,l)\in\Lam}||\zeta||_{2}^2-\sum_{(k,l)\in\Lam}\la
  (I-\ma N)^{-1}\zeta,\ma N\zeta\ra= \sum_{(k,l)\in\Lam}||\Pi_{2}(\frac{\chi_{1}^{k}\chi_{2}^l}{\alpha})||_{2}^2+R_{1}+R_{2}.$ 
  Par inégalité triangulaire on a directement $|R_{1}|\le ||\ma W||_{\Lam} (||\ma W||_{\Lam}+2||\ma V||_{\Lam})$. 
  Par ailleurs, en réutilisant l'égalité (\ref{fuchs}), on obtient la majoration 
  $|R_{2}|=|\sum_{\Lam}|\la (I-\ma N)^{-1}\zeta,\ma N \zeta\ra|\le \sum_{\Lam}\la \zeta,\ma N \zeta\ra|
  +C \sum_{\Lam}|| \ma N\zeta||_{2}^2$ où $C$ est une constante au moins égale à $||(I-\ma N)^{-1}||$.
  Or un calcul direct à partir de la définition de $\ma N$ et de l'expression de $\zeta$ sous la forme 
  $\zeta(q)=\Pi_{2}(\frac{\alpha}{\bar\alpha}\Pi_{1}^\perp(\frac{q}{\bar\alpha}))$ donne les majoration et égalité suivantes : $||\ma N\zeta||_{2}^2 \le ||\Pi_{1}(\frac{\alpha}{\bar\alpha}\zeta)||_{2}^2$ et 
  $\la\zeta,\ma N \zeta\ra=||\Pi_{1}(\frac{\alpha}{\bar\alpha}\zeta)||_{2}^2$, ce qui permet de conclure.
  \end{preuve}
  Le théorème \ref{otlitchno} reste une traduction intrinsèque de la proposition suivante. 
  \begin{proposition}\label{pretheo}
  Le symbole $f$ vérifie les hypothèses du théorème \ref{clef}, on note $\sum_{(u,v)\in\ma{C}^+}\beta_{u,v}\chi_{1}^{u}\chi_{2}^v$ le développement en série de Fourier de $\frac{1}{\alpha}$ dans la factorisation  $f=\alpha\bar{\alpha}$. On note $l_{i}(\lambda))$ la longueur du côté $c_{i}$ du triangle $\Lam_{\lambda}$ et $(a_{i},b_{i})$ le vecteur normal extérieur \emph{unitaire} au côté $c_{i}$.
  Alors

$\trace\left( T_{\Lam_{\lambda}}(f)^{-1}\right)=\left|\left|\frac{1}{f}\right|\right|_{1}|\Lam_{\lambda}|-\sum_{i=1}^3l_{i}(\lambda))\sum_{(u,v)\in\ma{C}^+}|a_{i}u+b_{i}v|\;|\beta_{u,v}|^2+o(\lambda)$
\end{proposition}
La preuve de la proposition \ref{pretheo} est la conséquence de la liste des lemmes et corollaires qui commencent à partir du lemme \ref{motive}.
Le lemme suivant montre que la proposition \ref{pretheo}, entra\^{i}ne le théorème \ref{otlitchno} après avoir  remarquer que $\sgn(a_{i}u+b_{i}v)=\sgn((-1)^{i})$ pour $i\in\{1,2,3\}$.
\begin{lem}\label{intrinsou}
On reprend les hypothèses et notations de la proposition \ref{pretheo} et on suppose que $f$ vérifie les hypothèses du théorème \ref{clef}. De plus,  on suppose que $\frac{1}{f}$ est essentiellement bornée sur $\tore^2$. Alors :
\begin{enumerate}
\item Les familles $\{u\;|\beta_{u,v}|^2\}_{\ma C^+}$ et $\{v\;|\beta_{u,v}|^2\}_{\ma C^+}$ sont sommables.
\item  On a les égalités :
\begin{align} 
\sum_{(u,v)\in\ma{C}^+}u|\beta_{u,v}|^2&=\frac{1}{2}\sum_{(u,v)\in\zz^2} |u|\;\;\widehat{\log\left(\frac{1}{ f}\right)}(u,v)\;\overline{\widehat{\left(\frac{1}{f}\right)}(u,v)}\label{intrin}\\
\sum_{(u,v)\in\ma{C}^+}v|\beta_{u,v}|^2&=\frac{1}{2}\sum_{(u,v)\in\zz^2} |v|\;\;\widehat{\log\left( \frac{1}{f}\right)}(u,v)\;\overline{\widehat{\left(\frac{1}{f}\right)}(u,v)}.\label{intrintrin}
\end{align}
 \end{enumerate}
\end{lem}
\begin{preuve}{du lemme \ref{intrinsou}}
Considérons la fonction $\beta(z)=P[\frac{1}{\alpha}](z)$, holomorphe sur $U^2$ (voir les notations de la section \ref{nonotes}). En posant pour $w\in\tore^2$,  $\beta_r(w)=P[\frac{1}{\alpha}](rw)$, on sait que $\lim_{r\to 1}||\beta_r(w)-\frac{1}{\alpha}||_2=0$ (voir \cite{RUD} page $18$). On en déduit le développement en série entière de $\beta(z)$ à savoir 
\begin{equation}\label{xiao}
\beta(z)=\sum_{\nn^2}\widehat{\left(\frac{1}{\alpha}\right)}(u,v)z_1^{u}z_2^v=\sum_{\nn^2}\beta_{uv}z_1^{u}z_2^v=\sum_{\ma C^{+}}\beta_{uv}z_1^{u}z_2^v,
\end{equation}
 La première égalité est la traduction du noyau de Poisson (voir \cite{RUD} page $17$). 
La deuxième égalité provenant de la définition des coefficients $\beta_{u,v}$.
La troisième égalité de (\ref{xiao}) découle du théorème \ref{clef}.
En  partant de l'équation( \ref{xiao}), on obtient alors l'égalité 
\begin{equation}\label{xiaopin}
 \widehat{\left(\frac{1}{\overline\alpha}\right)}(-u_{0},-v_{0})=\overline{\widehat{\left(\frac{1}{\alpha}\right)}}(u_{0},v_{0}).
 \end{equation}
Profitant de l'analycité de $\beta(z)$, on a $\frac{\partial \beta_{r}}{\partial \theta_{2}}=i\sum_{\ma C^+}v\beta_{u,v}r^{u+v}\chi_{1}^{u}\chi_{2}^v$, ce qui donne si l'on note $f_{r}^{-1}=\beta_{r}\bar\beta_{r}$ :
\begin{align}\label{chop1}
-\sum_{\ma C^+}v|\beta_{u,v}|^2r^{2(u+v)}&=-i\la\frac{\partial \beta_{r}}{\partial \theta_{2}},\beta_{r}\ra_{2}=-i\la 
\frac{1}{\beta_{r}}\frac{\partial \beta_{r}}{\partial \theta_{2}},\frac{1}{f_{r}}\ra_{2}=-i\la \frac{\partial \log \beta_{r}}{\partial \theta_{2}},\frac{1}{f_{r}}\ra_{2}=\nonumber\\
&=\sum_{\ma C^+}v\widehat{\log \beta_{r}}(u,v)\overline{\widehat{\left(\frac{1}{f_{r}}\right)}(u,v)}=\sum_{\zz^2}v\widehat{\log \beta_{r}}(u,v)\overline{\widehat{\left(\frac{1}{f_{r}}\right)}(u,v)}
\end{align}
en utilisant  le théorème de Parseval et le fait que $\widehat{\log \beta_{r}}(u,v)=0$ en dehors de $\ma C^+$. Notons que 
 La fonction $\frac{1}{f_{r}}$ étant réelle, on a :  $\widehat{\left(\frac{1}{f_{r}}\right)}(u,v)=\overline{\widehat{\left(\frac{1}{f_{r}}\right)}(-u,-v)}$. \`{A} partir de cette dernière  égalité et de l'égalité (\ref{xiaopin}), on obtient en conjugant les égalités (\ref{chop1}), 
 \begin{align}\label{chopine}
-\sum_{\ma C^+}v|\beta_{u,v}|^2r^{2(u+v)}&=i\overline{\la\frac{\partial \beta_{r}}{\partial \theta_{2}},\beta_{r}\ra_{2}}=i\la 
\frac{1}{\bar\beta_{r}}\frac{\partial\overline{ \beta_{r}}}{\partial \theta_{2}},\frac{1}{f_{r}}\ra_{2}=i\la \frac{\partial \log \overline{\beta_{r}}}{\partial \theta_{2}},\frac{1}{f_{r}}\ra_{2}=\nonumber\\
&=\sum_{\ma C^+}v\widehat{\log \overline{\beta_{r}}}(-u,-v)\overline{\widehat{\left(\frac{1}{f_{r}}\right)}(-u,-v)}=\sum_{-\ma C^+}|v|\;\widehat{\log \overline{\beta_{r}}}(u,v)\overline{\widehat{\left(\frac{1}{f_{r}}\right)}(u,v)}\nonumber\\
&=\sum_{\zz^2}v\widehat{\log \overline{\beta_{r}}}(u,v)\overline{\widehat{\left(\frac{1}{f_{r}}\right)}(u,v)}.
\end{align}

En additionnant les égalités (\ref{chop1}) et (\ref{chopine}) et en faisant tendre $r$ vers $1$,  on obtient l'égalité (\ref{intrintrin}). M\^{e}me démonstration pour (\ref{intrin}).
\end{preuve}
Le lemme \ref{motive} donne une estimation de la somme des deuxième et troisième termes de la formule d'inversion (\ref{trasa}) utilisant les projecteurs $p_{i}, i=1,2,3$ liés aux côtés du triangles.
  \begin{lem}\label{motive}
  Rappelons que pour tout $i\in\{1,2,3\}$, l'opérateur $p_{i}$ désigne la projection orthogonale de $L^2(\tore^2)$ sur $\ma P(S_{i}^-)$ (voir les notations de la section \ref{notar}). Alors 
  \begin{align}\label{devagar}
 & \sum_{(k,l)\in\Lam_{\lambda}}\left|\left|\Pi_{1}(\frac{\chi_{1}^{k}\chi_{2}^l}{\bar\alpha})\right|\right|_{2}^2+\sum_{(k,l)\in\Lam_{\lambda}}\left|\left|\Pi_{2}(\frac{\chi_{1}^{k}\chi_{2}^l}{\alpha})\right|\right|_{2}^2 =
 \sum_{(k,l)\in\Lam_{\lambda}}\left|\left|p_{1}(\frac{\chi_{1}^{k}\chi_{2}^l}{\bar\alpha})\right|\right|_{2}^2+\sum_{(k,l)\in\Lam_{\lambda}}\left|\left|p_{3}(\frac{\chi_{1}^{k}\chi_{2}^l}{\bar\alpha})\right|\right|_{2}^2\\&+\sum_{(k,l)\in\Lam_{\lambda}}\left|\left|p_{2}(\frac{\chi_{1}^{k}\chi_{2}^l}{\alpha})\right|\right|_{2}^2-\sum_{(k,l)\in\Lam_{\lambda}}\left|\left|p_{1}p_{3}(\frac{\chi_{1}^{k}\chi_{2}^l}{\bar\alpha})\right|\right|_{2}^2.\nonumber
 \end{align}
   \end{lem}
   \begin{preuve}{du lemme \ref{motive}}
   Si on pose $\Sigma_{1}=S_{1}^-\cap\{(u,v)\in\zz^2\; v>0\},\Sigma_{3}=S_{3}^-\cap\{(u,v)\in\ma \zz^2\; ; \;\alpha_{1}u+\beta_{1 }v<0\},\Sigma_{13}=\mathbb{S}_{1}^-\setminus (\Sigma_{1}\cup\Sigma_{3})$,
  on peut écrire $\mathbb{S}_{1}^-$ comme la réunion disjointe $\Sigma_{1}\cup\Sigma_{3}\cup\Sigma_{13}$.
   On note respectivement par $p_{11},p_{13}$ les projections orthogonales de $L^2(\tore^2)$ sur respectivement $\ma P(\Sigma_{1}),\ma P(\Sigma_{3})$. Alors pour tout vecteur $\xi$ de 
   $L^2(\tore^2)$, on a l'égalité :
    \begin{equation}\label{hop}
   ||p_{11}(\xi)||_{2}^2=||p_{1}(\xi)||_{2}^2-||p_{1}p_{3}(\xi)||_{2}^2
    \end{equation}
    En effet,
    
   $\la (p_{1}-p_{1}p_{3})(\xi),(p_{1}-p_{1}p_{3})(\xi)\ra=||p_{1}(\xi)||_{2}^2+||p_{1}p_{3}(\xi)||_{2}^2-\la p_{1}(\xi),
   p_{1}p_{3}(\xi)\ra-\la  p_{1}p_{3}(\xi),p_{1}(\xi)\ra$, les deux derniers termes étant égaux à $||p_{1}p_{3}(\xi)||_{2}^2$ car $p_{1}$ et $p_{3}$ commutent et on conclut. De m\^{e}me $||p_{13}(\xi)||_{2}^2=||p_{3}(\xi)||_{2}^2-||p_{1}p_{3}(\xi)||_{2}^2$. Par ailleurs la projection orthogonale sur $\ma P(\Sigma_{13})$ est $p_{1}p_{3}$.
   On a donc  
   \begin{align*}
 & \sum_{(k,l)\in\Lam_{\lambda}}\left|\left|\Pi_{1}(\frac{\chi_{1}^{k}\chi_{2}^l}{\bar\alpha})\right|\right|_{2}^2=  \sum_{(k,l)\in\Lam_{\lambda}}\left|\left|p_{11}(\frac{\chi_{1}^{k}\chi_{2}^l}{\bar\alpha})\right|\right|_{2}^2+ \sum_{(k,l)\in\Lam_{\lambda}}\left|\left|p_{13}(\frac{\chi_{1}^{k}\chi_{2}^l}{\bar\alpha})\right|\right|_{2}^2+ \sum_{(k,l)\in\Lam_{\lambda}}\left|\left|p_{1}p_{3}(\frac{\chi_{1}^{k}\chi_{2}^l}{\bar\alpha})\right|\right|_{2}^2\\
  & (\text{car les ensembles}\; \Sigma_{1},\Sigma_{3} \Sigma_{13}\; \text{sont disjoints})\\
 & =\sum_{(k,l)\in\Lam_{\lambda}}\left|\left|p_{1}(\frac{\chi_{1}^{k}\chi_{2}^l}{\bar\alpha})\right|\right|_{2}^2+ \sum_{(k,l)\in\Lam_{\lambda}}\left|\left|p_{3}(\frac{\chi_{1}^{k}\chi_{2}^l}{\bar\alpha})\right|\right|_{2}^2- \sum_{(k,l)\in\Lam_{\lambda}}\left|\left|p_{1}p_{3}(\frac{\chi_{1}^{k}\chi_{2}^l}{\bar\alpha})\right|\right|_{2}^2
  \end{align*}
  On conclut en notant que $\Pi_{2}=p_{2}$.
   \end{preuve}
Les lemmes \ref{fondecal} et \ref{frondemer} vont donner une estimation asymptotique des termes de droite de 
l'équation (\ref{devagar}).
\begin{lem}(lemme technique 1)\label{fondecal}\\
On pose $m=\lfloor \sqrt \lambda\rfloor$ ( c'est à dire la partie entière de $\sqrt\lambda$) et on considère les deux points $B_{1}$ et $B_{2}$ suivants  du triangle $\Lam_{\lambda}$ :  $B_{1}$ est le point du côté $c_{2}$ tel que les segments $OA_{1}$ et $B_{2}B_{1}$ soient parallèles et $B_{2}=(ma,0)$ (voir figure \ref{calculus}).

On note $\ma{D}_{1}$ le quadrilatère convexe $OA_{1}B_{1}B_{2}$ inclus dans $\Lam_{\lambda}$, $\ma B_{1}$ la bande incluse dans $ \ma C^+$ délimitée par les droites $(0A_{1})$ et $(B_{2}B_{1})$ (voir figure \ref{calculus}). 
Dans la factorisation $f=\alpha\bar\alpha$ du symbole $f$, on a $\frac{1}{\alpha}\in H^\ma C$ et on pose 
$\frac{1}{\alpha}=\sum_{\ma C}\beta_{u,v}\chi_{1}^{u}\chi_{2}^v$ son développement en série de Fourier. Alors
\begin{align}
&1)\sum_{(k,l)\in \ma{D}_{1}}\left|\left|p_{1}\left(\frac{\chi_{1}^{k}\chi_{2}^l}{\overline{\alpha}}\right)\right|\right|_{2}^2=\lambda \Big(\sum_{(u,v)\in\ma{B}_{1}}|\alpha_{1}u+\beta_{1}v||\beta_{u,v}|^2\Big)+o(\lambda)\label{flor1}\\
&2)\sum_{(k,l)\in \Lam_{\lambda}\setminus{\ma{D}_{1}}}\left|\left|p_{1}\left(\frac{\chi_{1}^{k}\chi_{2}^l}{\overline{\alpha}}\right)\right|\right|_{2}^2 =o(\lambda)\label{flor2}
\end{align}
\end{lem}
\begin{coro}{}\label{odecal}

Avec les notations du lemme \ref{fondecal}, on a l'estimation asymptotique 
\[\sum_{(k,l)\in \Lam_{\lambda}}\left|\left|p_{1}\left(\frac{\chi_{1}^{k}\chi_{2}^l}{\overline{\alpha}}\right)\right|\right|_{2}^2=\lambda \Big(\sum_{(u,v)\in\ma{B}_{1}}|\alpha_{1}u+\beta_{1}v||\beta_{u,v}|^2\Big)+o(\lambda)\]
\end{coro}
\begin{figure}
\unitlength 0.9mm
\begin{picture}(105, 85)
\thinlines
\put(20,0){\line(0,1){85}}
\put(0,20){\line(1,0){105}}
\put(57.5,20){\line(1,1){25}}%
 \put(57.5,20){\line(1,1){15}}
\thicklines
 \put(20,20){\line(1,0){75}}
 \put(20,20){\line(1,1){50}}
  \put(20,20){\line(1,1){50}}
 \put(70,70){\line(1,-2){25}} 
\dottedline{1}(35,35)( 72.5,35)

 \dottedline{1}(83,45)( 117,79)%
 \dottedline{1}(94,20)(146,72)
 \dottedline{1}(70,70)(145,70)
 \put(17,21){$0$}
 \put(58,16){$B_2$}
 \put(93,16){$A_2$}
 \put(63,69){$A_1$}
  \put(30,36){$A_0$}%
   \put(74,33){$B_3$}%
  \put(109,66){$C_1$}
 \put(86,46){$B_1$}
 \put(146,66){$C_2$}
 \put(59,45){$\mathcal D_{1}$}
 \put(76,25){$\Lambda\setminus\mathcal D_{1}$}
 \put(94,72){$\mathcal B_1$}
 \end{picture}
 \caption{lemmes techniques $1$ et $3$}
 \label{calculus}
 \end{figure}

\begin{preuve}{du lemme \ref{fondecal}}
\emph{Démonstration du premier item}.

Le côté $c_{1}$ est porté par la droite d'équation $\alpha_{1}x+\beta_{1}y=0$, avec $\alpha_{1}<0,\beta_{1}>0$. La droite parallèle à $c_{1}$ qui porte le segment $B_{1}B_{2}$ a pour équation $\alpha_{1}x+\beta_{1}y=ma \alpha_{1}$ . Chaque point de $\ma{D}_{1}$ se trouve sur une droite $D_{i}$ d'équation $\alpha_{1}x+\beta_{1}y=-i, \;0\le i\le ma|\alpha_{1}|.$ 
On a l'inclusion $\ma{D}_{1}\subset \{(k,l)\;;\;-m a\alpha_{1}<\alpha_{1}k+\beta_{1}l\le 0,\;l\le -\lambda\alpha_{1}\}$.
Si le point $(k,l)$ appartient à $ D_{i}$, alors
$(k-u,l-v)$ est dans $ S_{1}^-$ si et seulement si $\alpha_{1}u+\beta_{1}v<-i.$ \emph{Cette condition ne dépend que de $D_{i}$ et non de $(k,l)\in D_{i}$}. Posons  pour tout $i\in\nn\cap \big[0,ma|\alpha_{1}|\big]$
\[\Delta_{i}=\{(u,v)\in\ma{C}^+\;;\;\alpha_{1}u+\beta_{1}v<-i\}.\]
On a alors $\Delta_{i+1}\subset \Delta_{i}$ et  $\Delta_{i}=\bigcup_{j>i}D_{j}\cap \mathcal{C}^+$. Alors, pour tout $(k,l)\in D_{i}$, on a  $ p_{1}(\chi_{1}^{k-u}\chi_{2}^{l-v})\not= 0$ si et seulement si $(u,v)\in \Delta_{i}$. On a donc \begin{equation}\label{floflo}
\sum_{(k,l)\in \ma{D}_{1}}\left|\left|p_{1}\left(\myfrac{\chi_{1}^k\chi_{2}^l}{\overline{\alpha}_{}}\right)\right|\right|^2=\sum_{i=0}^{ma|\alpha_{1}|}\sum_{(k,l)\in D_{i}\cap\ma{D}_{1}}\;\;\;\sum_{(u,v)\in\Delta_{i}}|\beta_{uv}|^2.\end{equation}
Pour estimer le nombre de points de $D_{i}$, noté $|D_{i}|$, examinons deux cas, suivant que $i$ soit divisible par $\alpha_{1}$ ou non.
Pour tout entier $i$ divisible par $\alpha_{1}$ de $[0,ma\alpha_{1}]$, on a, compte tenu des coordonnées de $A_{1}$ et $A_{2}$ (voir sous-section \ref{notar}) l'encadrement  $|D_{ma\alpha_{1}}\cap\ma{D}_{1}|=(|D_{0}\cap\ma{D}_{1}|-m)+1\le |D_{i}\cap\ma{D}_{1}|\le|D_{0}\cap\ma{D}_{1}|=\lambda+1$. L'égalité $|D_{ma}\cap\ma{D}_{1}|=(|D_{0}\cap\ma{D}_{1}|-m)+1$ résulte d'un calcul élémentaire et l'encadrement de $|D_{i}|$  pour $i$ divisible par $\alpha_{1}$ vient du fait que les points entiers de $D_{i}$ sont de la forme $(x+k,y)$ si $i=k\alpha_{1}$ et $(x,y)\in D_{0}$. Supposons maintenant que  $i\in \big]s|\alpha_{1}|,(s+1)|\alpha_{1}|\big[$ et $0\le s< ma$ et 
notons $u$ et $v$ des coefficients de Bézout tels que $\alpha_{1}u+\beta_{1}v=1$. Remarquons  que ces coefficients sont positifs car $\alpha_{1}<0$. La droite $D_{i}\cap\zz^2$ est alors la famille $\{(x_{n},y_{n}\}_{n\in\zz}$ avec $x_{n}=n\beta_{1}-ui$ et $y_{n}=n|\alpha_{1}|-vi$. Comme la droite $D_{s|\alpha_{1}|}$ rencontre le côté $c_{2}$ en un point d'ordonnée $\lambda|\alpha_{1}|+\frac{s\alpha_{1}}{a}$, les points $(x_{n},y_{n})$ de $D_{i}\cap \Lam_{\lambda}$ vérifient la contrainte $0\le y_{n}\le  \lambda|\alpha_{1}|+\frac{s\alpha_{1}}{a}$. Ce qui donne l'encadrement de $n$ suivant : $\frac{vi}{|\alpha_{1}|}\le n\le \frac{vi}{|\alpha_{1}|}+\lambda-\frac{s}{a}$  qui montre compte tenu de la majoration de $s$ que $\lambda-m\le|D_{i}|\le\lambda$.
Ainsi, pour un choix de $m$ approprié (par exemple en prenant $m$ égal à la partie entière par défaut de $\lambda^{1/2}$) on a pour tout entier $i$ de $[0,ma\alpha_{1}]$ l'égalité  $|D_{i}|=(\lambda+1)+o(\lambda)$. Par conséquent en notant $o(1)$ un infiniment petit quand $\lambda$ tend vers $+\infty$, on a :
\begin{align*}
&\sum_{(k,l)\in \ma{D}_{1}}\left|\left|p_{1}\left(\myfrac{\chi_{1}^k\chi_{2}^l}{\overline{\alpha}}\right)\right|\right|^2=(\lambda+1)\sum_{i=0}^{ma|\alpha_{1}|}\sum_{(u,v)\in\Delta_{i}}|\beta_{uv}|^2+o(\lambda)\sum_{i=0}^{ma|\alpha_{1}|}\sum_{(u,v)\in\Delta_{i}}|\beta_{uv}|^2.
\end{align*}
\'{E}valuons le terme $\sum_{i=0}^{ma|\alpha_{1}|}\sum_{(u,v)\in\Delta_{i}}|\beta_{uv}|^2$.
\begin{align*}
&\sum_{i=0}^{ma|\alpha_{1}|}\sum_{(u,v)\in\Delta_{i}}|\beta_{uv}|^2=\sum_{\Delta_{0}}|\beta_{uv}|^2+\sum_{\Delta_{1}}|\beta_{uv}|^2+\ldots+\sum_{\Delta_{ma|\alpha_{1}|}}|\beta_{uv}|^2\\
&=(ma|\alpha_{1}|+1)\sum_{\Delta_{ma|\alpha_{1}|}\cap\mathcal{C}^+}|\beta_{uv}|^2+(ma|\alpha_{1}|)\sum_{D_{ma|\alpha_{1}|}\cap\mathcal{C}^+}|\beta_{uv}|^2+\ldots+\sum_{D_{1}\cap\mathcal{C}^{+}}|\beta_{uv}|^2
\end{align*}
On en déduit

$$\sum_{(k,l)\in \ma{D}_{1}}\left|\left|p_{1}\left(\myfrac{\chi_{1}^k\chi_{2}^l}{\overline{\alpha}_{}}\right)\right|\right|^2=(\lambda+1)\left((ma|\alpha_{1}|+1)\sum_{\Delta_{ma|\alpha_{1}|}\cap\mathcal{C}^+}|\beta_{uv}|^2+\sum_{i=1}^{ma|\alpha_{1}|}\sum_{(u,v)\in D_{i}}|i||\beta_{u,v}|^2+o(1)\right).$$
 Or si $(u,v)\in\Delta_{ma|\alpha_{1}|}\cap\mathcal{C}^+$, on a 
$u> ma$,  donc 

$(ma|\alpha_{1}|+1)\sum_{\Delta_{ma|\alpha_{1}|}\cap\mathcal{C}^+}|\beta_{uv}|^2\le
(|\alpha_{1}|+\frac{1}{ma}) \sum_{\Delta_{ma|\alpha_{1}|}\cap\mathcal{C}^+}u|\beta_{uv}|^2=o(1)$,
 ce qui permet de conclure directement le premier item.

 \noindent \emph{Démontstration du deuxième item}.
 
 En reprenant la m\^{e}me démarche que pour le premier item on aboutit à l'égalité suivante analogue à l'égalité (\ref{floflo}) :
 \begin{equation*}\label{}
\sum_{(k,l)\in \Lam_{\lambda}\setminus{\ma{D}_{1}}}\left|\left|p_{1}\left(\myfrac{\chi_{1}^k\chi_{2}^l}{\overline{\alpha}_{}}\right)\right|\right|^2=\sum_{i=ma|\alpha_{1}|+1}^{\lambda a}\sum_{(k,l)\in D_{i}}\sum_{(u,v)\in\Delta(D_{i})}|\beta_{uv}|^2,\end{equation*}
Or pour $i>ma|\alpha_{1}|$ on a $\left|D_{i}\cap(\Lam_{\lambda}\setminus\ma{D}_{1})\right|
\le (\lambda-m)|\alpha_{1}|<\lambda|\alpha_{1}|$ ce qui conduit, avec la m\^{e}me calcul que celui du premier point, à l'inégalité suivante :

$\sum_{(k,l)\in \Lam_{\lambda}\setminus{\ma{D}_{1}}}\left|\left|p_{1}\left(\myfrac{\chi_{1}^k\chi_{2}^l}{\overline{\alpha}_{}}\right)\right|\right|^2\le ( \lambda-m)|\alpha_{1}|\sum _{(u,v)\in\mathcal C^+\setminus\ma B_{1}}|\alpha_{1}u+\beta_{1}v||\beta_{u,v}|^2.$
Et on conclut en notant que $\sum _{(u,v)\in\mathcal C^+\setminus\ma B_{1}}|\alpha_{1}u+\beta_{1}v||\beta_{u,v}|^2=o(1)$ et $\lambda-m=O(\lambda)$.
\end{preuve}

\begin{figure}
\unitlength 0.9mm
\begin{picture}(105, 85)
\thinlines
\put(20,0){\line(0,1){85}}
\put(0,20){\line(1,0){105}}
\thicklines
  
  \put(30,20){\line(1,1){70}}
   \put(74,20){\line(1,1){65}}%
   
 \put(20,70){\line(1,0){130}} %
\put(20,45){\line(1,0){120}} 

 \put(17,21){$0$}
 \put(30,16){$ma$}
\put(74,16){$na$}
  \put(9,45){$p|\alpha_{1}|$}
   \put(9,70){$q|\alpha_{1}|$}
  \put(96,75){$\Lam_{(m,n,q,\infty)}$}
 \put(80,60){$\Lam_{(m,n,p,q)}$}
  \put(125,60){$\Lam_{(n,\infty,p,q)}$}
 \end{picture}
 \caption{lemme technique $2$}\label{lemtec}
 \end{figure}
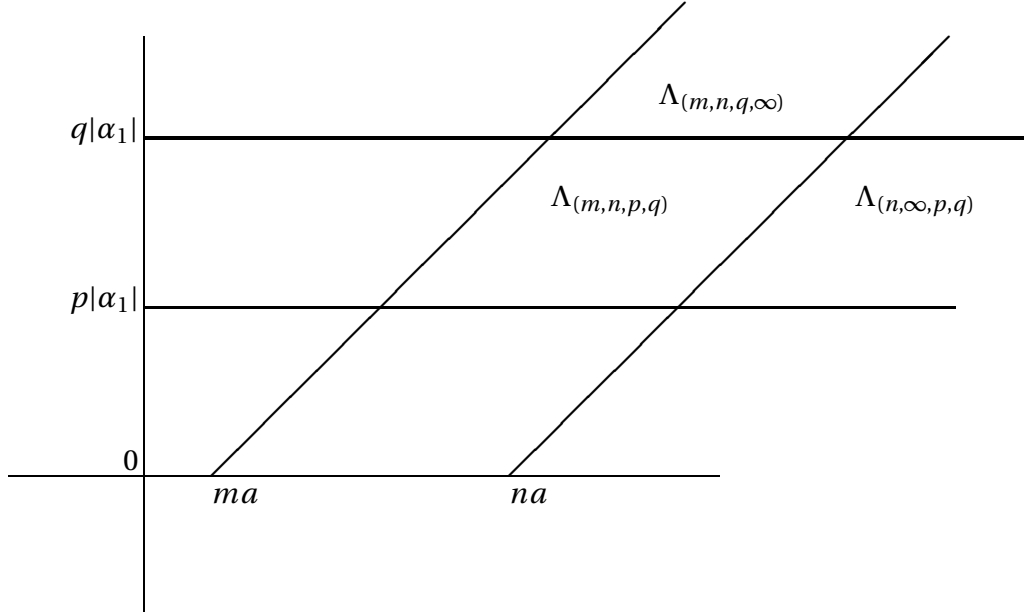

\begin{lem}{(lemme technique 2)}\label{frondemer}

Conformément aux notations précédentes, $D_{i}$ désigne la droite d'équation $\alpha_{1}x+\beta_{1}y=-i$ où $i$ est un entier naturel et $\Delta_{i} =\bigcup_{j>i}D_{j}$. On se donne le quadruplet $(m,n,p,q)\in\nn^4$ tel que 
$0\le m<n\le\lambda$ et $0\le p<q\le\lambda$. Il définit le parallélogramme $\Lam_{(m,n,p,q)} $ délimité par les droites $D_{ma\alpha_{1}},D_{na\alpha_{1}}$ et les droites horizontales d'équations respectives $y=p|\alpha_{1}|$ et $y=q|\alpha_{1}|$, illustré par la figure \ref{lemtec}.
On posera de plus  :
\begin{enumerate}
\item $D_{i}(v_{0})=\{(u,v)\in D_{i}\;;\;v\ge v_{0}\}$,
\item  $\Delta_{i}(v_{0})=\bigcup_{j>i}D_{j}(v_{0})$,
\item $ \Lam_{(m,n,q,\infty)}=\ma C^+\bigcap\left(\bigcup_{j=ma|\alpha_{1}|+1}^{na|\alpha_{1}|}D_{j}(q|\alpha_{1}|)\right)$ ( voir figure \ref{lemtec}),
\item $\Lam_{(n,\infty,pq)}= \ma C^+\bigcap\left\{ (u,v)\in\nn^2\;;\; v\in\big[p|\alpha_{1}|+1,q|\alpha_{1}|\big]\right\}\bigcap \Delta_{na|\alpha_{1}|}$ ( voir figure \ref{lemtec}).
\end{enumerate}
Alors $$\sum_{(k,l)\in\Lam_{(m,n,p,q)}}\left|\left|p_{1}p_{3}\left(\frac{\chi_{1}^{k}\chi_{2}^l}{\overline{\alpha}}\right)\right|\right|_{2}^2=\sum_{i=1}^5T_{i}$$ où les termes $T_{i} $ sont des fonctions de $(m,n,p,q)$ définies par:
\begin{align*}
T_{1}&=\big((n-m)|\alpha_{1}|+1\big)\big((q-p)|\alpha_{1}|+1\big)\sum_{\Delta_{na|\alpha_{1}|}(q|\alpha_{1}|)}|\beta_{u,v}|^2, \\
T_{2}&=\big((q-p)|\alpha_{1}|+1\big)\sum_{ \Lam_{(m,n,q,\infty)}}\big(|\alpha_{1}u+\beta_{1}v|-ma|\alpha_{1}|\big)|\beta_{u,v}|^2,\\
T_{3}&=\big((n-m)|\alpha_{1}|+1\big)\sum_{ \Lam_{(n,\infty,p,q)}}\big(v-p|\alpha_{1}|\big)|\beta_{u,v}|^2,\\
T_{4}&=\sum_{(u,v)\in \Lam_{(m,n,p,q)}}|\alpha_{1}u+\beta_{1}v|(v-p|\alpha_{1}|)|\beta_{u,v}|^2\\
T_{5}&=-ma|\alpha_{1}|\sum_{(u,v)\in \Lam_{(m,n,p,q)}}(v-p|\alpha_{1}|)|\beta_{u,v}|^2.
\end{align*}
\end{lem}

\begin{coro}\label{fond}
On a l'évaluation asymptotique suivante quand $\lambda$ tend vers $+\infty$ :
\[\sum_{\Lam_{\lambda}}\left|\left|p_{1}p_{3}\left(\frac{\chi_{1}^{k}\chi_{2}^l}{\overline{\alpha}}\right)\right|\right|_{2}^2=o(\lambda).\]
\end{coro}
\begin{preuve}{ du lemme \ref{frondemer}.}
Tout point $(k,l)\in\Lam_{(m,n,p,q)}$ appartient à une droite $D_{i}$ où $i$ est un entier de $\big[ ma\alpha_{1}, na\alpha_{1}\big]$. Pour un point $(k,l)$ de $D_{i}$ et tout $(u,v)\in \ma C$, on remarque que $(k-u,l-v)\in S_{1,0}^-\cup S_{3,0}^-=S_{1}^-\cup S_{3}^-$ si et seulement si $(u,v)\in\Delta_{i}(l)$.

Cette remarque permet de décrire la somme $\sum_{(k,l)\in\Lam_{(m,n,p,q)}}\left|\left|p_{1}p_{3}\left(\frac{\chi_{1}^{k}\chi_{2}^l}{\overline{\alpha}}\right)\right|\right|_{2}^2$ sous la forme :
\[
\sum_{(k,l)\in\Lam_{(m,n,p,q)}}\left|\left|p_{1}p_{3}\left(\frac{\chi_{1}^{k}\chi_{2}^l}{\overline{\alpha}}\right)\right|\right|_{2}^2=
\sum_{i=ma|\alpha_{1}|}^{na|\alpha_{1}|}\;\;\sum_{(k,l)\in D_{i}\cap\Lam_{(m,n,p,q)}}\;\;\sum_{(u,v)\in\Delta_{i}(l)}|\beta_{u,v}|^2.\]
Posons $\Gamma_{i}=\sum_{(k,l)D_{i}\cap\Lam_{(m,n,p,q)}}\;\;\sum_{(u,v)\in\Delta_{i}(l)}|\beta_{u,v}|^2.$

En remarquant que $\Delta_{j}(l+1)\subset \Delta_{j}(l)$, on a  par un calcul analogue à celui de l'item $1$ du lemme \ref{fondecal} :
$\Gamma_{i}=\sum_{\Delta_{i}(p|\alpha_{1}|})|\beta_{u,v}|^2+\ldots+\sum_{\Delta_{i}(q|\alpha_{1}|)}|\beta_{u,v}|
^2= \big((q-p)|\alpha_{1}|+1\big)\sum_{\Delta_{i}(q|\alpha_{1}|)}|\beta_{u,v}|^2
+\sum_{(u,v)\in \Delta_{i}}
\sum_{v=p|\alpha_{1}\lambda+1}^{q|\alpha_{1}|}(v-p|\alpha_{1}|)|\beta_{u,v}|^2$. On en déduit,
en posant :
 \[\tau_1=\big((q-p)|\alpha_{1}|+1\big)\sum_{i=ma|\alpha_1|}^{na|\alpha_1|}\sum_{\Delta_{i}(q|\alpha_{1}|)}|\beta_{u,v}|^2\] et
  \[\tau_2=\sum_{i=ma|\alpha_1|}^{na|\alpha_1|}\sum_{(u,v)\in \Delta_{i}}
\sum_{v=p|\alpha_{1}|\lambda+1}^{q|\alpha_{1}|}(v-p|\alpha_{1}|)|\beta_{u,v}|^2,\]
l'égalité suivante : 
$\sum_{(k,l)\in\Lam_{(m,n,p,q)}}\left|\left|p_{1}p_{3}\left(\frac{\chi_{1}^{k}\chi_{2}^l}{\overline{\alpha}}\right)\right|\right|_{2}^2=\tau_1 +\tau_2$.
 
\noindent\emph{\'{E}valuation de} $\frac{1}{\big((q-p)|\alpha_{1}|+1\big)}\tau_1.$\\

On a : $\frac{1}{\big((q-p)|\alpha_{1}|+1\big)}\tau_1=\sum_{_{\Delta_{ma|\alpha_1|}}\big(q|\alpha_1|\big)}|\beta_{u,v}|^2+\ldots +\sum_{_{\Delta_{na|\alpha_1|}}(q|\alpha_1|)}|\beta_{u,v}|^2$

$=\big((n-m)a|\alpha_1|+1\big)\sum_{_{\Delta_{na|\alpha_1|}}(q|\alpha_1|)}|\beta_{u,v}|^2+\displaystyle\sum_{j=ma|\alpha_1|+1}^{na|\alpha_1|}\;\;\sum_{(u,v)\in D_{_j}(q|\alpha_1|)}(j-ma|\alpha_{1}|)|\beta_{u,v}|^2$

$=\big((n-m)a|\alpha_1|+1\big)\sum_{_{\Delta_{na|\alpha_1|}}(q|\alpha_1|)}|\beta_{u,v}|^2+
\sum_{\Lam_{(m,n,q,\infty)}}\big(|\alpha_1 u+\beta_1 v|-ma|\alpha_1|\big)\big|\beta_{u,v}\big|^2$, par définition de $D_{_j}(q|\alpha_1|)$ et de $\Lam_{(m,n,q,\infty)}$.
On a donc $\tau_{1}=T_{1}+T_{2}$.

\noindent\emph{\'{E}valuation de} $\tau_2.$\\
En posant $\gamma_v=\sum_{i=ma|\alpha_1|}^{na|\alpha_1|}\;\;\sum_{(u,v)\in\Delta_i}(v-p|\alpha_1|)|\beta_{u,v}|^2$ on peut écrire 
$\tau_2=\sum_{v=p|\alpha_1|+1}^{q|\alpha_1|}\gamma_v$.
La remarque utilisée déjà dans le lemme \ref{fondecal}, à savoir $\Delta_{i+1}\subset\Delta_i$ permet d'obtenir une estimation de $\gamma_v$ par un calcul analogue à celui  de $\Gamma_{i}$ :

\noindent$\gamma_v=\sum_{(u,v)\in\Delta_{ma|\alpha_1|}} (v-p|\alpha_1|)+\ldots+\sum_{(u,v)\in\Delta_{na|\alpha_1|}}(v-p|\alpha_1|)=$

\noindent$=\big((n-m)a|\alpha_1|+1\big)\displaystyle\sum_{(u,v)\in\Delta_{na|\alpha_{1}|}}(v-p|\alpha_1|)|\beta_{u,v}|^2+
\displaystyle\sum_{j=ma|\alpha_{1}|}^{na|\alpha_{1}|}\displaystyle\sum_{(u,v)\in D_{j}}(j-ma|\alpha_{1}|)(v-p|\alpha_{1}|)\beta_{u,v}|^2$.
Ainsi, en posant 

$\nu_{1}=\big((n-m)a|\alpha_1|+1\big)\displaystyle\sum_{v=p|\alpha_{1}|+1}^{q|\alpha_{1}|}\;\displaystyle\sum_{(u,v)\in\Delta_{na|\alpha_{1}|}}(v-p|\alpha_1|)|\beta_{u,v}|^2$ et

$\nu_{2}=\displaystyle\sum_{v=p|\alpha_{1}|+1}^{q|\alpha_{1}|}\;\displaystyle\sum_{j=ma|\alpha_{1}|}^{na|\alpha_{1}|}\;\displaystyle\sum_{(u,v)\in D_{j}}(j-ma|\alpha_{1}|)(v-p|\alpha_{1}|)|\beta_{u,v}|^2$, on obtient  $\tau_{2}=\nu_{1}+\nu_{2}$.
Reste à remarquer que $\nu_{1}=T_{3}$ et $\nu_{2}=T_{4}+T_{5}$ et que $\tau_{2}=T_{3}+T_{4}+T_{5}$.
Le compte y est.
\end{preuve}
\begin{preuve}{ du corollaire \ref{fond}.}
Soit $A_{0}$ le point du côté $OA_{1}$ d'ordonnée $m\alpha_{1}$ et $B_{3}$ le point intérieur au segment $B_{2}B_{1}$ d'ordonnée $m\alpha_{1}$ (voir figure \ref{calculus}).
Le polygone $\Lam_{\lambda}$ est inclus dans la réunion  des trois parallélogrammes $OA_{0}B_{3}B_{2}$, $A_{0}A_{1}C_{1}B_{3}$ et 
$B_{2}C_{1}C_{2}A_{2}$ (voir figure \ref{calculus}) où, rappelons le, $B_{2}$ est le point de coordonnées $(ma,0)$, $0\le m\le \lambda$. Par définition de $\Lam_{(m,n,p,q)}$, on a : $OA_{0}B_{3}B_{2}=\Lam_{(0,m,0,m)}$
$A_{0}A_{1}C_{1}B_{3}=\Lam_{(0,m,m,\lambda)}$ et $B_{2}C_{1}C_{2}A_{2}=\Lam_{(m,\lambda,0,\lambda)}$
et on a la majoration 
 $\sum_{(k,l)\in\Lam_{\lambda}}\left|\left|p_{1}p_{3}\left(\frac{\chi_{1}^{k}\chi_{2}^l}{\overline{\alpha}}\right)\right|\right|_{2}^2$
 
 \noindent$\le\sum_{(k,l)\in\Lam_{(0,m,m,\lambda)}}\left|\left|p_{1}p_{3}\left(\frac{\chi_{1}^{k}\chi_{2}^l}{\overline{\alpha}}\right)\right|\right|_{2}^2+\sum_{(k,l)\in\Lam_{(m,\lambda,0,\lambda})}\left|\left|p_{1}p_{3}\left(\frac{\chi_{1}^{k}\chi_{2}^l}{\overline{\alpha}}\right)\right|\right|_{2}^2+\sum_{(k,l)\in\Lam_{(0,m,0,m)}}\left|\left|p_{1}p_{3}\left(\frac{\chi_{1}^{k}\chi_{2}^l}{\overline{\alpha}}\right)\right|\right|_{2}^2.$
 Chacune des deux sommes de droite relève du lemme \ref{frondemer}. En choisissant $m=\lfloor\sqrt{\lambda}\rfloor$, on obtient le tableau suivant qui résume le comportement asymptotique des termes $\{T_{i}\}_{i=1,\ldots,5}$ définis dans le lemme \ref{frondemer} et en fin de compte des trois termes du majorant de $\sum_{(k,l)\in\Lam_{\lambda}}\left|\left|p_{1}p_{3}\left(\frac{\chi_{1}^{k}\chi_{2}^l}{\overline{\alpha}}\right)\right|\right|_{2}^2$  lorsque $\lambda$ tend vers l'infini.
\[
\begin{array}{|c||c|c|c|}
\hline (m,n,p,q)&(0,\lfloor\sqrt{\lambda}\rfloor,\lfloor\sqrt{\lambda}\rfloor,\lambda)&(\lfloor\sqrt{\lambda}\rfloor,\lambda,0,\lambda)&(0,\lfloor\sqrt{\lambda}\rfloor,0,\lfloor\sqrt{\lambda}\rfloor)\\
\hline \hline T_{1}&\sqrt{\lambda}\;o(1)&\lambda\;o(1)&\sqrt{\lambda}\;o(1)\\
\hline T_{2}&\lambda \;o(1)&\lambda \;o(1)&\sqrt{\lambda}\;o(1)\\
\hline T_{3}&{\lambda}\;o(1)&\lambda \;o(1)&\sqrt{\lambda}\;o(1)\\
\hline T_{4}&\lambda\;o(1)&\lambda \;o(1)&\lambda\;O(1)\\
\hline T_{5}&0&\sqrt{\lambda}\;o(1)&0\\
\hline \sum_{i=1}^5T_{i}&o(\lambda)&o(\lambda)&o(\lambda)\\
\hline
\end{array}
\]
Montrons que $\sum_{(k,l)\in\Lam_{(0,m,m,\lambda)}}\left|\left|p_{1}p_{3}\left(\frac{\chi_{1}^{k}\chi_{2}^l}{\overline{\alpha}}\right)\right|\right|_{2}^2=o(\lambda)$ (deuxième colonne du tableau).

On a :

1)
\begin{align*}
T_{1}(0,\lfloor\sqrt{\lambda}\rfloor,\lfloor\sqrt{\lambda}\rfloor,\lambda)&=(\lfloor\sqrt{\lambda}\rfloor|\alpha_{1}|+1)\big((\lambda-\sqrt{\lambda})\big|\alpha_{1}|+1\big)\sum_{\Delta_{_{\lfloor\sqrt{\lambda}\rfloor a|\alpha_{1}|}}(\lambda |\alpha_{1}|)}|\beta_{u,v}|^2\\
&\le (\lfloor\sqrt{\lambda}\rfloor|\alpha_{1}|+1)\sum_{\Delta_{_{\lfloor\sqrt{\lambda}\rfloor a|\alpha_{1}|}}(\lambda |\alpha_{1}|)}v|\beta_{u,v}|^2
\end{align*}
puisque $v>\lambda |\alpha_{1}|$ lorsque $(u,v)\in \Delta_{_{\lfloor\sqrt{\lambda}\rfloor a|\alpha_{1}|}}(\lambda|\alpha_{1}|)$. Compte tenu du lemme \ref{intrinsou} et des hypothèses du théorème \ref{otlitchno}  , la somme $\sum_{\Delta_{_{\lfloor\sqrt{\lambda}\rfloor a|\alpha_{1}|}}(\lambda |\alpha_{1}|)}v|\beta_{u,v}|^2$ est le reste d'une série convergente et on conclut : $T_{1}(0,\lfloor\sqrt{\lambda}\rfloor,0,\lambda)=\sqrt{\lambda}\;o(1)$.

2) $T_{2}(0,\lfloor\sqrt{\lambda}\rfloor,\lfloor\sqrt{\lambda}\rfloor,\lambda)=\big((\lambda-\sqrt{\lambda})|\alpha_{1}|+1\big)\sum_{ \Lam_{(0,\lfloor\sqrt{\lambda}\rfloor,\lambda,\infty)}}\big|\alpha_{1}u+\beta_{1}v\big|\big|\beta_{u,v}|^2={\lambda}\;o(1)$, car la somme est comme ci-dessus un reste de série convergente.

3) Pour la m\^{e}me raison, $T_{3}(0,\lfloor\sqrt{\lambda}\rfloor,\lfloor\sqrt{\lambda}\rfloor,\lambda)=\big(\lfloor\sqrt{\lambda}\rfloor|\alpha_{1}|+1\big)\sum_{ \Lam_{(\lfloor\sqrt{\lambda}\rfloor,\infty,\lfloor\sqrt{\lambda}\rfloor,\lambda)}}\big(v-\lfloor\sqrt{\lambda}\rfloor|\alpha_{1}|\big)|\beta_{u,v}|^2$.

$={\lambda}\;o(1)$

4)En remarquant que $\frac{\big(v-\lfloor\sqrt{\lambda}\rfloor|\alpha_{1}|\big)}{\lambda}=O(1)$ sur 
$\Lam_{(0,\lfloor\sqrt{\lambda}\rfloor,\lfloor\sqrt{\lambda}\rfloor,\lambda)}$, on a 
$T_{4}(0,\lfloor\sqrt{\lambda}\rfloor,\lfloor\sqrt{\lambda}\rfloor,\lambda)=\sum_{(u,v)\in \Lam_{(0,\lfloor\sqrt{\lambda}\rfloor,\lfloor\sqrt{\lambda}\rfloor,\lambda)}}|\alpha_{1}u+\beta_{1}v|(v-\lfloor\sqrt{\lambda}\rfloor|\alpha_{1}|)|\beta_{u,v}|^2=\lambda O(1) o(1)=\lambda o(1)$, ce qui justifie les éléments de la deuxième colonne, puisque $T_{5}(0,\lfloor\sqrt{\lambda}\rfloor,\lfloor\sqrt{\lambda}\rfloor,\lambda)=0$.
La justification des termes des  troisième et quatrième colonnes est du m\^{e}me ordre.
\end{preuve}

Nous allons enfin énoncer un lemme qui aboutira à l'estimation des restes $R_{1}$ et $R_{2}$ qui apparaissent dans le développement de la trace défini par l'égalité (\ref{trasa}) du corollaire \ref{tras}.
\begin{lem}{(Lemme technique\; $3$)}\label{resttec}
En plus des notations du corollaire \ref{tras}, on pose $Z$ l'opérateur défini sur $\ma P(\Lam_{\lambda})$ par 
$Z(q)=\Pi_{1}\left(\frac{\alpha}{\bar\alpha}\Pi_{2}\big(\frac{q}{\alpha}\big)\right)$.
\begin{enumerate}
\item $||\ma W||_{\Lam_{\lambda}}=o(\sqrt\lambda)$.
\item $||Z||_{\Lam_{\lambda}}=o(\sqrt\lambda)$.
\end{enumerate}
\end{lem}
\begin{preuve}{du lemme \ref{resttec}}
\emph{Premier item}

On considère les projections $p_{11},p_{13}$ définies au début  de la preuve du lemme \ref{motive}.
En décomposant $\Pi_{1}$ sous la forme $\Pi_{1}=p_{11}+p_{13}+p_{1}p_{3}$,  (somme de trois projections orthogonales sur trois 
espace de polynômes de spectres deux à deux disjoints), on peut écrire par inégalité triangulaire sur les norme de Hilbert-Schmidt :
$$||\ma W||_{\Lam}\le||\ma W_{1}||_{\Lam}+||\ma W_{2}||_{\Lam}+||\ma W_{3}||_{\Lam}$$ 
avec
$\left\{
\begin{array}{cc}
  ||\ma W_{1}||^2_{\Lam}=\sum_{(k,l)\in\Lam}||\Pi_{2}\big(\frac{\bar\alpha}{\alpha}p_{11}(\frac{\chi_{1}^k\chi_{2}^l}{\bar\alpha}\big)||^2,&||\ma W_{2}||^2_{\Lam}=\sum_{(k,l)\in\Lam}||\Pi_{2}\big(\frac{\bar\alpha}{\alpha}p_{13}(\frac{\chi_{1}^k\chi_{2}^l}{\bar\alpha}\big)||^2 \\||\ma W_{3}||^2_{\Lam}=\sum_{(k,l)\in\Lam}||\Pi_{2}\big(\frac{\bar\alpha}{\alpha}p_{1}p_{3}(\frac{\chi_{1}^k\chi_{2}^l}{\bar\alpha}\big)||^2.&
  \end{array}
\right.$

\'{E}valuons $ ||\ma W_{1}||^2_{\Lam}$. On peut écrire conformément à la figure \ref{calculus}, $||\ma W_{1}||^2_{\Lam}= ||\ma W_{1}||^2_{\ma D_{1}}+ ||\ma W_{1}||^2_{\Lam\setminus\ma D_{1}}$.
Soit $\ma C_{m}$ le paraléllogramme $OA_{0}B_{3}B_{2}$ de la figure \ref{calculus} où $m$ est défini comme dans le lemme \ref{fondecal}.
On considère alors   le polynôme $\frac{1}{\alpha_{m}}=\sum_{(u,v)\in\ma C_{m}}\beta_{u,v}\chi_{1}^{u}\chi_{2}^v$.
Par l'inégalité triangulaire dans la norme $||.||_{\ma D_{1}}$, on a 
$||\ma W||_{\ma D_{1}}\le\left(\sum_{\ma D_{1}}\left|\left|\Pi_{2}\big(\frac{\bar\alpha}{\alpha_{m}}p_{11}(\frac{\chi_{1}^k\chi_{2}^l}{\bar\alpha})\big)\right|\right|_{2}^2\right)^{1/2}+\left(\sum_{\ma D_{1}}\left|\left|\Pi_{2}\Big(\bar\alpha\big(\frac{1}{\alpha}-\frac{1}{\alpha_{m}}\big)p_{11}(\frac{\chi_{1}^k\chi_{2}^l}{\bar\alpha})\Big)\right|\right|_{2}^2\right)^{1/2}$. Or d'une part, on a sur $\ma D_{1}$ l'égalité
$\Pi_{2}\big(\frac{\bar\alpha}{\alpha_{m}}p_{11}(\frac{\chi_{1}^k\chi_{2}^l}{\bar\alpha})\big)=0$ et d'autre part la majoration

$\left(\sum_{\ma D_{1}}\left|\left|\Pi_{2}\Big(\bar\alpha\big(\frac{1}{\alpha}-\frac{1}{\alpha_{m}}\big)p_{11}(\frac{\chi_{1}^k\chi_{2}^l}{\bar\alpha})\Big)\right|\right|_{2}^2\right)^{1/2}\le$

$\le ||\alpha||_{\infty}||\frac{1}{\alpha}-\frac{1}{\alpha_{m}} ||_{\infty}\left(\sum_{\ma D_{1}}\left|\left|p_{11}(\frac{\chi_{1}^k\chi_{2}^l}{\bar\alpha})\right|\right|_{2}^2\right)^{1/2}
\le
 ||\alpha||_{\infty}||\frac{1}{\alpha}-\frac{1}{\alpha_{m}} ||_{\infty}\left(\sum_{\ma D_{1}}\left|\left|p_{1}(\frac{\chi_{1}^k\chi_{2}^l}{\bar\alpha})\right|\right|_{2}^2\right)^{1/2}$ grâce à l'égalité (\ref{hop}) de la démonstration du lemme \ref{motive}.
 On déduit du premier point du lemme \ref{fondecal} que $||\ma W||_{\ma D_{1}}=o(1)O(\sqrt\lambda)=o(\sqrt\lambda)$.
 Par ailleurs $||\ma W_{1}||^2_{\Lam\setminus\ma D_{1}}\le \sum_{\Lam\setminus\ma D_{1}}||p_{1}(\frac{\chi_{1}^k\chi_{2}^l}{\bar\alpha})||^2=o(\lambda)$ d'après le deuxième point du lemme \ref{fondecal}. On a en fin de compte $ ||\ma W_{1}||_{\Lam}=o({\sqrt\lambda})$.
 Pour des raisons de symétries on a également $ ||\ma W_{2}||_{\Lam}=o(\sqrt\lambda).$
 Reste à évaluer $||\ma W_{3}||_{\Lam}$. Par le  corollaire \ref{fond}, on a :

 $||\ma W_{3}||^2_{\Lam}= \displaystyle\sum_{(k,l)\in\Lam}\left|\left|\Pi_{2}\big(\frac{\bar\alpha}{\alpha}p_{1}p_{3}(\frac{\chi_{1}^k\chi_{2}^l}{\bar\alpha}\big)\right|\right|_{2}^2\le\sum_{\Lam}\left|\left|p_{1}p_{3}(\frac{\chi_{1}^k\chi_{2}^l}{\bar\alpha})\right|\right|_{2}^2=o(\lambda)$ et on conclut.

 \emph{Deuxième item.}
 \begin{figure}
\unitlength 0.9mm
\begin{picture}(105, 85)
\thinlines
\put(20,0){\line(0,1){85}}
\put(0,20){\line(1,0){105}}
\thicklines
 \put(20,20){\line(1,0){75}}
 \put(20,20){\line(1,1){50}}
  \put(20,20){\line(1,1){50}}
 \put(70,70){\line(1,-2){25}} 
 \put(64,64){\line(1,-2){22}} 
 %
 \dottedline{1}(28,28)( 90,28)
 \dottedline{1}(30,20)( 73,63)%
 \put(17,21){$0$}
 \put(30,16){$B_2$}
 \put(93,16){$A_2$}
 \put(66,72){$A_1$}
  \put(82,16){$A_4$}
   \put(93,29){$A_5$}   \put(76,29){$A_6$}
   \put(57,64){$A_3$}
 \put(22,29){$A_0$} \put(42,29){$B_{4}$}
   \put(65,50){$B_3$}
    \put(65,62){$D_{12}$}%
     \put(85,24){$D_{13}$}%
      \put(75,46){$D_{2}$}%
 \put(75,64){$B_1$}
 \put(55,38){$\mathcal D_{2}$}
  \put(32,43){$c_{1}$} \put(90,43){$c_{2}$} \put(57,15){$c_{3}$}
 \put(26,23){$\mathcal {C}_m$}
 \end{picture}
 \caption{calcul 2}\label{calcululu}
 \end{figure}

 Afin d'obtenir des majorations suffisamment fines, on considère la partition  de $\Lam_{\lambda}$  représentée par la figure \ref{calcululu} où $A_{0}=(m|\alpha_{1}|,m|\alpha_{1}|),A_{5}$ est le point de $c_{2}$ d'ordonnée $m|\alpha_{1}|$, $B_{2}=(ma,0),A_{4}=((\lambda-m)a,0), A_{2}=(\lambda,0)$, le segment 
 $A_{4}A_{3}$ est parallèle au côté $c_{2}$, le segment 
 $B_{2}B_{1}$ est parallèle au côté $c_{1}$ . Enfin $B_{4}$ est l'intersection des droites $(A_{0}A_{5})$ et
 $(B_{1}B_{2})$. 
 
 Alors $D_{12}$ est l'intersection avec $\zz^2$ du parallélogramme de $\rr^2$ fermé $A_{3}A_{1}B_{1}B_{3}$, de m\^{e}me $D_{13}$ est l'intersection avec $\zz^2$  du parallélogramme fermé $A_{2}A_{4}A_{6}A_{5}$, on note $D_{2}$  l'intersection avec $\zz^2$ du quadrilatère fermé $B_{3}B_{1}A_{5}A_{6}$. On a alors
 $\ma D_{2}= \Lam_{\lambda}\setminus (D_{12}\cup D_{2}\cup D_{13})$. Pour finir $\ma C_{m}$ est l'intersection avec $\zz^2$ du parallélogramme de $\rr^2$ fermé $OA_{0}B_{4}B_{2}$.
 
 \'{E}crivons la décomposons $||Z||^2_{\Lam_{\lambda}}$ sous la forme : 
 $||Z||^2_{\Lam_{\lambda}}=||Z||^2_{D_{12}}+||Z||^2_{D_{13}}+||Z||^2_{D_{2}}+||Z||^2_{\ma D_{2}}$. 
 La démarche pour l'évaluation des termes $||Z||^2_{D_{12}}, ||Z||^2_{D_{13}},||Z||^2_{\ma D_{2}}$ est calquée en partie sur le lemme \ref{fondecal}. Cette démarche ne s'appliquera pas pour le terme $||Z||^2_{ D_{2}}$ car elle conduit alors à une estimation trop grossière. Nous développons ce calcul dans les deux points suivants.
 
\noindent 1) \emph{Estimation de} $||Z||^2_{D_{12}}, ||Z||^2_{D_{13}},||Z||^2_{\ma D_{2}}$.

On a pour ces trois quantités les majorations : 

\noindent$||Z||^2_{\ma D_{2}}\le \sum_{(k,l)\in\ma D_{2}}\left|\left|\Pi_{2}\big(\frac{\chi_{1}^k\chi_{2}^l}{\alpha}\big)\right|\right|_{2}^2,||Z||^2_{D_{12}}\le \sum_{(k,l)\in D_{12}}\left|\left|\Pi_{2}\big(\frac{\chi_{1}^k\chi_{2}^l}{\alpha}\big)\right|\right|_{2}^2,||Z||^2_{D_{13}}\le \sum_{(k,l)\in D_{13}}\left|\left|\Pi_{2}\big(\frac{\chi_{1}^k\chi_{2}^l}{\alpha}\big)\right|\right|_{2}^2$.

En vue d'estimer les trois majorants  précédents, introduisons quelques notations. 

Le côté $c_{2}$ est porté, conformément à la figure \ref{triangle} par la droite d'équation $\alpha_{2}x+\beta_{2}y=\lambda a\alpha_{2}$ et le segment$A_{3}A_{4}$ par la droite  $\alpha_{2}x+\beta_{2}y=(\lambda-m) a\alpha_{2}$. Tout point entier $(k,l)$ de 
 $\Lam_{\lambda}$ appartient à une droite $d_{i}$ d'équation $\alpha_{2}x+\beta_{2}y=i$, où $i$ est un entier naturel vérifiant :
\[
 \begin{cases}
 (\lambda-m)a\alpha_{2}\le i\le\lambda a\alpha_{2}&\text{ si}\;d_{i}\cap( \Lam_{\lambda}\setminus \ma D_{2})\not=\emptyset \\
 0\le i\le (\lambda-m)a\alpha_{2}&\text{si}\; d_{i}\cap\ma D_{2}\not=\emptyset.
 \end{cases}
\]
 Si on pose $\Gamma_{i}=\left\{(u,v)\in \ma C^+\;;\; \alpha_{2}u+\beta_{2}v>\lambda a\alpha_{2}-i\right\},$ alors un point entier $(k,l)$ de $d_{i}$ vérifie $\Pi_{2}(\chi_{1}^k\chi_{2}^l)\not= 0$ si et seulement si $(u,v)\in\Gamma_{i}$. 
On peut écrire maintenant de fa\c con analogue au début de la démonstration du lemme \ref{fondecal} :

$\sum_{(k,l)\in\ma D_{2}}\left|\left|\Pi_{2}\big(\frac{\chi_{1}^k\chi_{2}^l}{\alpha}\big)\right|\right|_{2}^2=\sum_{i=0}^{(\lambda-m)a\alpha_{2}}\sum_{(k,l)\in d_{i}\cap\ma D_{2}}\sum_{(u,v)\in\Gamma_{i}}|\beta_{u,v}|^2.$
Or $|d_{i}\cap\ma D_{2}|\le O(\lambda-m)$, d'où $\sum_{(k,l)\in\ma D_{2}}\left|\left|\Pi_{2}\big(\frac{\chi_{1}^k\chi_{2}^l}{\alpha}\big)\right|\right|_{2}^2\le O(\lambda-m)\sum_{i=0}^{(\lambda-m)a\alpha_{2}}\sum_{(u,v)\in\Gamma_{i}}|\beta_{u,v}|^2$. Et par le m\^{e}me procédé de sommation que dans la fin de la démonstration du premier item du lemme
\ref{fondecal}, mais en tenant compte maintenant des inclusions $\Gamma_{i}\subset \Gamma_{i+1}$, on obtient : $\sum_{(k,l)\in\ma D_{2}}\left|\left|\Pi_{2}\big(\frac{\chi_{1}^k\chi_{2}^l}{\alpha}\big)\right|\right|_{2}^2\le $

$\le O(\lambda-m)\big((\lambda-m)a\alpha_{2}+1\big)\sum_{\Gamma_{0}}|\beta_{u,v}|^2+O(\lambda-m)\sum_{i=ma\alpha_{2}+1}^{\lambda a\alpha_{2}}(i-ma\alpha_{2})\sum_{(u,v)\in d_{i}}|\beta_{u,v}|^2$.
Par ailleurs on a la majoration :

 $0\le\sum_{i=ma\alpha_{2}+1}^{\lambda a\alpha_{2}}(i-ma\alpha_{2})\sum_{(u,v)\in d_{i}}|\beta_{u,v}|^2\le
\sum_{i=ma\alpha_{2}+1}^{\lambda a\alpha_{2}}\sum_{(u,v)\in d_{i}}i|\beta_{u,v}|^2=$

$=\sum_{i=ma\alpha_{2}+1}^{\lambda a\alpha_{2}}\sum_{(u,v)\in d_{i}}(\alpha_{2}u+\beta_{2}v)|\beta_{u,v}|^2$. En fin de compte $\sum_{(k,l)\in\ma D_{2}}||\Pi_{2}\big(\frac{\chi_{1}^k\chi_{2}^l}{\alpha}\big)||^2\le t_{1}+t_{2}$ où $t_{1}=O(\lambda-m)\big((\lambda-m)a\alpha_{2}+1\big)\sum_{\Gamma_{0}}|\beta_{u,v}|^2$ et 
$t_{2}=O(\lambda-m)\sum_{\ma D_{m,\lambda}}(\alpha_{2}u+\beta_{2}v)|\beta_{u,v}|^2$ si on désigne par 
$\ma D_{m,\lambda}$ la bande incluse dans $\Lam_{\lambda}$ entre les droites $d_{ma\alpha_{2}+1}$ et $d_{\lambda a\alpha_{2}}$.

Donnons maintenant une estimation asymptotique de $t_{1} $ et $t_{2}$ en prenant $m=\lfloor \sqrt{\lambda}\rfloor$.

On a l'inclusion $\Gamma_{0}\subset \mathbb S_{2}^- $ et par conséquent  si $(u,v)\in\Gamma_{0}$, on a 
$\alpha_{2}u+\beta_{2}v\ge \lambda a\alpha_{2}$. Ainsi $t_{1}\le O(\lambda-m)\sum_{\ma C^+\cap \mathbb S_{2}^-}(\alpha_{2}u+\beta_{2}v)|\beta_{u,v}|^2=O(\lambda-m)o(1)=o(\lambda)$ par le premier item du lemme \ref{intrinsou}.
De la m\^{e}me manière $t_{2}=o(\lambda)$ par le choix de $m$ et le premier item lemme \ref{intrinsou}.

\noindent Ainsi $\sum_{(k,l)\in\ma D_{2}}\left|\left|\Pi_{2}\big(\frac{\chi_{1}^k\chi_{2}^l}{\alpha}\big)\right|\right|_{2}^2=o(\lambda).$
 
\noindent Avec la m\^{e}me idée, on obtient pour la somme $\sum_{(k,l)\in D_{12}}\left|\left|\Pi_{2}\big(\frac{\chi_{1}^k\chi_{2}^l}{\alpha}\big)\right|\right|_{2}^2$ :

\noindent$\sum_{(k,l)\in D_{12}}\left|\left|\Pi_{2}\big(\frac{\chi_{1}^k\chi_{2}^l}{\alpha}\big)\right|\right|_{2}^2=\sum_{i=(\lambda-m)a\alpha_{2}}^{\lambda a\alpha_{2}}\sum_{(k,l)\in d_{i}\cap D_{12}}\sum_{(u,v)\in\Gamma_{i}}|\beta_{u,v}|^2 $ et la majoration $|d_{i}\cap D_{12}|\le O(m)$ conduit par un calcul parallèle au calcul précédent à l'estimation : $\sum_{(k,l)\in D_{12}}\left|\left|\Pi_{2}\big(\frac{\chi_{1}^k\chi_{2}^l}{\alpha}\big)\right|\right|_{2}^2 $

$\le O(m)\;ma\alpha_{2}\sum_{(u,v)\in\Gamma_{(\lambda-m)a\alpha_{2}}}|\beta_{u,v}|^2+O(m)\;\sum_{i=1}^{m a\alpha_{2}}\sum_{(u,v)\in d_{i}}(\alpha_{2}u+\beta_{2}v)|\beta_{u,v}|^2$. Or si $(u,v)\in\Gamma_{(\lambda-m)a\alpha_{2}}$, on a $\alpha_{2}u+\beta_{2}v>ma\alpha_{2}$. Par conséquent on peut écrire :

\noindent$O(m)\;ma\alpha_{2}\sum_{(u,v)\in\Gamma_{(\lambda-m)a\alpha_{2}}}|\beta_{u,v}|^2
<O(m)\sum_{(u,v)\in\Gamma_{(\lambda-m)a\alpha_{2}}}(\alpha_{2}u+\beta_{2}v)|\beta_{u,v}|^2=o(1)$, avec toujours le choix $m=\lfloor \sqrt{\lambda}\rfloor$ et le premier item du lemme \ref{intrinsou}. De plus, $\sum_{(u,v)\in d_{i}}(\alpha_{2}u+\beta_{2}v)|\beta_{u,v}|^2 =O(1)$, toujours d'après le premier item du lemme \ref{intrinsou}. 
On en conclut que \begin{equation}\label{ouf}
\sum_{(k,l)\in D_{12}}\left|\left|\Pi_{2}\big(\frac{\chi_{1}^k\chi_{2}^l}{\alpha}\big)\right|\right|_{2}^2=O(\lfloor \sqrt{\lambda}\rfloor)=o(\lambda).
\end{equation} 
Pour des raisons de symétrie on a également 
$\sum_{(k,l)\in D_{13}}\left|\left|\Pi_{2}\big(\frac{\chi_{1}^k\chi_{2}^l}{\alpha}\big)\right|\right|_{2}^2=o(\lambda).$

\noindent 2) \emph{Estimation de} $||Z||^2_{D_{2}}$.

On a $||Z||^2_{D_{2}}=\sum_{D_{2}}||\Pi_{1}\left(\frac{\alpha}{\bar\alpha}\Pi_{2}\big(\frac{q}{\alpha}\big)\right)||^2$.
On pose $\frac{1}{\alpha_{m}}=\displaystyle\sum_{(u,v)\in\ma C_{m}}\beta_{u,v}\chi_{1}^k\chi_{2}^l$. On écrit alors :

$\sum_{D_{2}}||\Pi_{1}\left(\frac{\alpha}{\bar\alpha}\Pi_{2}\big(\frac{q}{\alpha}\big)\right)||^2=
\sum_{D_{2}}||\Pi_{1}\left(\frac{\alpha}{\bar\alpha_{m}}\Pi_{2}\big(\frac{q}{\alpha}\big)\right)+
\Pi_{1}\left(\alpha\big(\frac{1}{\bar\alpha}-\frac{1}{\bar\alpha_{m}}\big)\Pi_{2}\big(\frac{q}{\alpha}\big)||^2\right)$.

\noindent Or pour tout $(k,l)\in D_{2}$, on a $\Pi_{1}\left(\frac{\alpha}{\bar\alpha_{m}}\Pi_{2}\big(\frac{q}{\alpha}\big)\right)=0$ et par conséquent 

\noindent$\sum_{D_{2}}||\Pi_{1}\left(\frac{\alpha}{\bar\alpha}\Pi_{2}\big(\frac{q}{\alpha}\big)\right)||^2\le ||\alpha||_{\infty}\left|\left|\frac{1}{\bar\alpha}-\frac{1}{\bar\alpha_{m}}\right|\right|_{\infty}\sum_{D_{2}}||\Pi_{2}(\frac{\chi_{1}^k\chi_{2}^l}{\alpha})||^2=o(1)O(\lambda-2m)$, d'après un calcul analogue à celui conduisant à l'égalité (\ref{ouf}). On en conclut que  $\sum_{D_{2}}||\Pi_{1}\left(\frac{\alpha}{\bar\alpha}\Pi_{2}\big(\frac{q}{\alpha}\big)\right)||^2=o(\lambda)$.
L'item $2$ est  démontré.
\end{preuve}
\begin{coro}\label{last}
Avec les notations du corollaire \ref{tras}, on a $R_{1}+R_{2}=o(\lambda)$.
\end{coro}
\begin{preuve}{du corollaire \ref{last}}
Le calcul de $||Z||_{D_{12}}$ qui est un calcul semblable à celui conduisant à l'égalité (\ref{flor1}), montre que 
$||\ma V||_{\Lam}=O(\sqrt{\lambda})$. On déduit donc du premier item du lemme \ref{resttec} que 
$|R_{1}|\le ||\ma W||_{\Lam}(||\ma W||_{\Lam}+2||\ma V||_{\Lam})=o(\lambda)$ 

Pour le reste $R_{2}$, remarquons que l'on peut écrire 

$\sum_{(k,l)\in\Lam_{\lambda}}\left|\left|\Pi_{1}\Big(\frac{\alpha}{\bar\alpha}
\zeta(\chi_{1}^k\chi_{2}^l)\Big)\right|\right|_{2}^2 =\sum_{(k,l)\in\Lam_{\lambda}}\left|\left|\Pi_{1}\left(\frac{\alpha}{\bar\alpha}\Pi_{2}
\big(\frac{\chi_{1}^k\chi_{2}^l}{\alpha}\big)\right)-\Pi_{1}\left(\frac{\alpha}{\bar\alpha}
\Pi_{2}\big(\frac{\bar\alpha}{\alpha}\Pi_{1}\big(\frac{\chi_{1}^k\chi_{2}^l}{\bar\alpha}\big)\big)\right)\right|\right|_{2}^2$
en rempla\c cant dans l'expression de $\zeta$ la projection $\Pi_{1}^\perp$ par $I-\Pi_{1}$. 
Mais alors par l'inégalité triangulaire sur les normes de Hilbert-Schmidt, on a la majoration :

$|R_{2}|^{1/2}\le ||Z||_{\Lam}+||\ma W||_{\Lam}$, en notant que $\left|\left|\Pi_{1}\left(\frac{\alpha}{\bar\alpha}
\Pi_{2}\big(\frac{\bar\alpha}{\alpha}\Pi_{1}\big(\frac{\chi_{1}^k\chi_{2}^l}{\bar\alpha}\big)\big)\right)\right|\right|_{2}^2\le 
\left|\left|\Pi_{2}\big(\frac{\bar\alpha}{\alpha}\Pi_{1}\big(\frac{\chi_{1}^k\chi_{2}^l}{\bar\alpha}\big)\big)\right|\right|_{2}^2$, 
ce qui permet de conclure avec le lemme \ref{resttec}.
\end{preuve}
\begin{preuve}{de la proposition \ref{pretheo}}
Elle découle directement et dans l'ordre indiqué du corollaire \ref{tras}, du lemme \ref{motive}, du corollaire \ref{odecal} et enfin du corollaire \ref{last}.
\end{preuve}
\section{Théorème du déterminant. Cas du triangle}\label{deter}
Rappelons l'énoncé du théorème du déterminant (théorème \ref{determi}).

\noindent\textbf{Théorème du déterminant}

\noindent On suppose, en plus des hypothèses du théorème de trace que $f\in\mathfrak K$ et $\ln f\in\mathfrak K$ (voir au début de la section \ref{deter}).

Posons 
$\mu_{1}(f)=-\frac{1}{2}\sum_{(u,v)\in\zz^2}|u| |\widehat{\ln f}(u,v)|^2$,
$  \mu_{2}(f)=-\frac{1}{2}\sum_{(u,v)\in\zz^2}|v| |\widehat{\ln f}(u,v)|^2$.
 Alors \[\det\;T_{\Lam_{\lambda}}(f)=e^{|\Lam_{\lambda}|\;||\ln f||_{1}}e^{-\lambda\big(\mathfrak S_{1}\mu_{1}(f)+\mathfrak S_{2}\mu_{2}(f)+o(1)\big)}.\]

\subsection{Preuve du théorème du déterminant}
On pose $f=1-h$ et  $f_{t}=1-th$ pour $t\in]0,1[$ avec $||h||_{\infty}<1.$

\noindent Preuve en quatre points.

\emph{Premier point.}  On \'{e}tablit l'égalité suivante :
 \begin{equation}\label{iopl}
\ln \left(\det \Big(T_{\Lam_{\lambda}}(1-h)\Big)\right)-\trace\Big(T_{\Lam_{\lambda}}(\ln(1-h))\Big)
=-\int_{0}^1\frac{1}{t}\trace\left(\Big(T_{\Lam_{\lambda}}(1-th)\Big)^{-1}-T_{\Lam_{\lambda}}\Big(\frac{1}{1-th}\Big)\right)dt.\end{equation}
Cela résulte des deux égalités suivantes :
\begin{align}
&\frac{d}{dt}\ln \left(\det \big(T_{\Lam_{\lambda}}(1-th)\big)\right)=\trace\left(\frac{1}{t}\Big(I-\big(T_{\Lam_{\lambda}}(1-th)\big)^\text{-1}\Big)\right)\label{lop}\\
&\frac{d}{dt}\trace\left(T_{\Lambda_{\lambda}}\big(\ln(1-th)\big)\right)=-\trace \left(T_{\Lambda_{\lambda}}\big(\frac{h}{1-th}\big)\right).\label{pol}
\end{align}

Pour prouver l'égalité (\ref{lop}), on montre tout d'abord, par une vérification directe l'identité suivante :
$\frac{1}{t}\left(I-\big(T_{\Lambda_{\lambda}}(1-th)\big)^{-1}\right)
=-T_{\Lambda_{\lambda}}(h)\big(T_{\Lambda_{\lambda}}(1-th)\big)^{-1}$.
Ensuite $T_{\Lambda_{\lambda}}(h)$ étant un opérateur hermitien, il est diagonalisable et on note $\{\nu_{i}\}_{i=1,\ldots,n_0}$ ses valeurs propres. Mais alors $\ln \left(\det \big(T_{\Lam_{\lambda}}(1-th)\big)\right)=
\sum_{i=0}^{n_{0}}\ln (1-t\nu_{i})$ d'où :

 $\frac{d}{dt}\ln \left(\det \big(T_{\Lam_{\lambda}}(1-th)\big)\right)=\frac{1}{t}\sum_{i=0}^{n_{0}}\left(1-\frac{1}{1-t\nu_{i}}\right)=\trace\left(\frac{1}{t}\Big(I-\big(T_{\Lam_{\lambda}}(1-th)\big)^\text{-1}\Big)\right).$
 L'égalité (\ref{pol}) est obtenue directement par le calcul de la dérivée partant du taux d'accroissement.
 En intégrant la différence des équations (\ref{lop}) et (\ref{pol}) et en remarquant que $I-T_{\Lam_{\lambda}}\big(\frac{th}{1-th}\big)=T_{\Lam_{\lambda}}\big(\frac{1}{1-th}\big)$, on obtient l'égalité (\ref{iopl}).
 
 \emph{Deuxième point}: on y met en évidence que le théorème est démontré si l'on établit l'équation (\ref{bisbis}).

 Commen\c cons par reformuler le théorème \ref{otlitchno}. Notons que $\sum_{i=1}^3l_{i}(\lambda)a_{i}
 =\lambda\sum_{i=1}^3a_{i}l_{i}(1)$, les nombres $l_{i}(1)$ représentant les longueurs des côtés du triangle de base correspondant à la valeur $1$ du paramètre. On pose alors :
 
 \begin{align}[left=\empheqlbrace]\label{couic}
&\forall t\in[0,1]\; f_{t}=1-th,\\
&c_{1}(f_{t})=\sum_{(u,v)\in\zz^2}|u|\;\widehat{\ln \frac{1}{f_{t}}}(u,v)\overline{\widehat{\frac{1}{f_{t}}}(u,v)},&c_{2}(f_{t})=\sum_{(u,v)\in\zz^2}|v|\;\widehat{\ln \frac{1}{f_{t}}}(u,v)\overline{\widehat{\frac{1}{f_{t}}}(u,v)}.\nonumber
\end{align}

 Le théorème de trace peut se réécrire ainsi pour le symbole $f_{t}$ :
 \begin{equation}\label{thrace}
 \trace\left(\big(T_{\Lam_{\lambda}}(f_{t})\big)^{-1}-T_{\Lam_{\lambda}}(\frac{1}{f_{t}})\right)=\lambda\big(\mathfrak S_{1}\;c_{1}(f_{t})+\mathfrak S_{2}\;c_{2}(f_{t})\big)+o(\lambda).
 \end{equation}
 Le théorème est démontré si l'on établit l'égalité suivante 
 \begin{equation}\label{bisbis}
 \mathfrak{S}_{1}\mu_{1}(f)+ \mathfrak{S}_{2}\mu_{2}(f)=\int_{0}^1\frac{ \mathfrak{S}_{1}c_{1}(f_{t})+\mathfrak{S}_{2}c_{2}(f_{t})}{t}dt
 \end{equation}
 En effet, des égalités (\ref{thrace}) et (\ref{bisbis}), on a directement 

$$\int_{0}^1\lim_{\lambda\to \infty}\frac{1}{\lambda} \trace\left(\big(T_{\Lam_{\lambda}}(f_{t})\big)^{-1}-T_{\Lam_{\lambda}}(\frac{1}{f_{t}})\right)\frac{dt}{t}=\mathfrak S_{1}\;\mu_{1}(f)+\mathfrak S_{2}\;\mu_{2}(f).$$
De la majoration pour $t$ assez petit et $\lambda$ assez grand, $\big|\frac{1}{\lambda} \trace\left(\big(T_{\Lam_{\lambda}}(f_{t})\big)^{-1}-T_{\Lam_{\lambda}}(\frac{1}{f_{t}})\right)\big|\le  2|c_{1}(f)+c_{2}(f)|$ déduite de (\ref{thrace}), on a par la convergence dominée 
\[\lim_{\lambda\to\infty}\int_{0}^1\frac{1}{\lambda} \trace\left(\big(T_{\Lam_{\lambda}}(f_{t})\big)^{-1}-T_{\Lam_{\lambda}}(\frac{1}{f_{t}})\right)\frac{dt}{t}=\mathfrak S_{1}\;\mu_{1}(f)+\mathfrak S_{2}\;\mu_{2}(f)\]
Ainsi l'égalité (\ref{iopl}) permet d'écrire :
$-\frac{1}{\lambda}\left(\ln \left(\det \Big(T_{\Lam_{\lambda}}(f)\Big)\right)-\trace\Big(T_{\Lam_{\lambda}}(f)\Big)\right)=\mathfrak S_{1}\;\mu_{1}(f)+\mathfrak S_{2}\;\mu_{2}(f)+o(1)$. Et on conclut.

Nous établirons l'égalité (\ref{bisbis}) dans le quatrième point de la démonstration. Elle nécessite deux lemmes 
que l'on énonce et démontre dans le troisième point.

\emph{Troisième point}: on établit quelques lemmes en vue de la démonstration de l'équation (\ref{bisbis}).

\begin{lem}\label{banane}
Avec les notations de la section \ref{trisym}, il existe sur $\ma L^\infty(\tore^2)$ une norme $N$ équivalente à $||.||_{\mathfrak K,\Lam}$, telle que la norme $||.||_{\ma B}$ définie pour tout $g\in L^\infty(\tore^2)$ par $||g||_{\ma B}=||g||_{\infty}+N(g)$ confère à $L^\infty(\tore^2)$ une structure d'algèbre de Banach.
\end{lem}
\begin{preuve}{du lemme \ref{banane}}\'{E}tablissons la preuve en deux étapes.

\emph{Première étape : existence de} $N$

Posons pour toute fonction $f\in L^2(\tore^2)$ : $N(f)=\left(\int_{\rr^2\times [0,2\pi]^2}\frac{|f(\theta+x)-f(\theta)|^2}{|| x||^3}dx d\theta\right)^{1/2}$ où
$\theta=(\theta_{1},\theta_{2})\in[0,\pi]^2,x=(x_{1},x_{2})\in\rr^2$ et où $f(\theta)$ est identifiée à $\tilde{f}(\theta)=f(e^{i\theta_{1}},e^{i\theta_{1}}).$ $N$ définit une norme sur $\ma L^\infty(\tore^2)$ ( voir définition de $\ma L^\infty(\tore^2)$ en sous-section \ref{trisym}). En appliquant l'égalité de Parceval à $\theta\mapsto f(\theta)$ et $\theta\mapsto f(\theta+x)$ on obtient directement :
\begin{equation}\label{thefirst}
N^2(f)=\sum_{\zz^2}|\widehat{f}(m,n)|^2\int_{\rr^2}\frac{|e^{i(mx_{1}+nx_{2})}-1|^2}{||x||^3}dx
\end{equation}
Notons $I(m,n)$  l'intégrale du membre de droite de l'égalité (\ref{thefirst}). On peut la réécrire sous la forme 
$I(m,n)= \int_{\rr^2}\frac{\sin^2(\frac{1}{2}(mx_{1}+nx_{2})}{(x_{1}^2+x_{2}^2)^{3/2}}dx_{1}dx_{2}$ qui met en évidence l'intégrabilité de la fonction $x\mapsto\frac{|e^{i(mx_{1}+nx_{2})}-1|^2}{||x||^3}$.
On a : \begin{equation}\label{joie}
 I(m,n)=K\;(m^2+n^2)^{1/2}
\end{equation} où $K$ est une  constante. En effet, le changement de variables en coordonnées polaires  $x_{1}=r\cos\theta, x_{2}=r\sin \theta$ donne 

\noindent$I(m,n)=\int_{[0,2\pi]\times [0,\infty[}\frac{\sin^2\left(\frac{1}{2}r(m\cos\theta+n\sin\theta)\right)}{r^2}dr d\theta=
\int_{0}^{2\pi}\big(\underbrace{\int_{0}^{+\infty}\frac{\sin^2\left(\frac{1}{2}r(m\cos\theta+n\sin\theta)\right)}{r^2}dr}_{I_{1}}\big)d\theta.$
Notons $\sin \alpha=\frac{m}{(m^2+n^2)^{1/2}}, \cos \alpha=\frac{n}{(m^2+n^2)^{1/2}}$ avec $0<\alpha<\pi$ et posons alors
$u=\frac{1}{2}r(m^2+n^2)^{1/2}$. Avec ce changement de variables, $I(m,n)$ s'écrit sous la forme indiquée par l'équation (\ref{joie}).

 Comme il existe deux constantes $C_{1}, C_{2}$ telles que 
 \[C_{1}(\mathfrak{S}_{1}|m|+\mathfrak{S}_{2}|n|)\le(m^2+n^2)^{1/2}\le C_{2}(\mathfrak{S}_{1}|m|+\mathfrak{S}_{2}|n|),\]
  les égalité (\ref{thefirst}) et (\ref{joie}) donnent l'équivalence des normes $N$ et $||.||_{\mathfrak{K},\Lam}$ sur 
 $\ma L^\infty(\tore^2)$.
 
 \emph{Deuxième étape : }$||.||_{\ma B}$ \emph{est une algèbre de Banach}
 
 Supposons démontrée l'inégalité suivante :
 \begin{equation}\label{algi}
 \forall f_{1},f_{2}\in L^\infty(\tore^2)\;N^2(f_{1}f_{2})=\big(||f_{1}||_{\infty}N(f_{2})+||f_{2}||_{\infty}N(f_{1})\big)^2
 \end{equation}
 Alors $||f_{1}f_{2}||_{\infty}+N(f_{1}f_{2})\le ||f_{1}||_{\infty}||f_{2}||_{\infty}+||f_{1}||_{\infty}N(f_{2})+||f_{2}||_{\infty}N(f_{1})\le \big(||f_{1}||_{\infty}+N(f_{1})\big)$
 
 \noindent$\big(||f_{2}||_{\infty}+N(f_{2})\big)=||f_{1}||_{\ma B}||f_{2}||_{\ma B}$, ce qui achève la deuxième étape.
 
 Reste à démontrer l'inégalité (\ref{algi}). 
 On a :
 
  $|(f_{1}f_{2})(\theta+x)-(f_{1}f_{2})(\theta)|^2=|f_{1}(\theta+x)\big(f_{2}(\theta+x)-f_{2}(\theta)\big)+
 f_{2}(\theta)\big(f_{1}(\theta+x)-f_{1}(\theta)\big)|^2$, dont on déduit l'inégalité 
 $|(f_{1}f_{2})(\theta+x)-(f_{1}f_{2})(\theta)|^2\le $
 
 \noindent$||f_{1}||^2_{\infty}|f_{2}(\theta+x)-f_{2}(\theta)|^2+
 ||f_{2}||^2_{\infty}|f_{1}(\theta+x)-f_{1}(\theta)|^2+2||f_{1}||_{\infty}||f_{2}||_{\infty}|f_{1}(\theta+x)-f_{1}(\theta)|\;|f_{2}(\theta+x)-f_{2}(\theta)|$, et par suite l'inégalité 
 
 \noindent$N^2(f_{1}f_{2})\le ||f_{1}||^2_{\infty}N^2(f_{2})+||f_{2}||^2_{\infty}N^2(f_{1})+2||f_{1}||_{\infty}||f_{2}||_{\infty}
 \int_{\rr^2\times \tore^2}\left|\frac{f_{1}(\theta+x)-f_{1}(\theta)}{||x||^{3/2}}\right|\left|\frac{f_{2}(\theta+x)-f_{2}(\theta)}{||x||^{3/2}}\right|dxd\theta$ et la démonstration s'achève avec l'inégalité de Hölder.
\end{preuve}
\begin{coro}{}\label{baba}

Il existe une constante $D$ telle que pour tout entier naturel $s$, on a $||h^s||_{\mathfrak K,\Lam}< D$. On peut choisir $D=1$.
\end{coro}
\begin{preuve}{ du corollaire \ref{baba}}
D'après le lemme \ref{banane}, il existe une constante $C$ telle que $N\le C||.||_{\mathfrak K,\Lam}$.
Quitte à diviser $f$ par une constante, on peut supposer, compte tenu de l'hypothèses $f\in \mathfrak K $ que $||h||_{\infty}+C ||h||_{\mathfrak K,\Lam}<1$. Mais alors pour un entier positif $s$, on a 
$||h^s||_{\mathfrak K,\Lam}\le ||h^s||_{\ma B}\le||h||^s_{\ma B}\le(||h||_{\infty}+C||h||_{\ma B})^s<1$.
\end{preuve}
\begin{lem}\label{ananas}
Pour $t\in]0,1[$ posons $\varphi_{t}=\ln (1-th)$. Alors $\lim_{\varepsilon\to 0}\big|\big|\frac{\varphi_{t+\varepsilon}-\varphi_{t}}{\varepsilon}-\frac{d\varphi_{t}}{dt}\big|\big|_{\mathfrak K,\Lam}=0$.
\end{lem}
\begin{preuve}{du lemme \ref{ananas}}
 En effet, $\lim_{\varepsilon\to 0}\big|\big|\frac{\varphi_{t+\varepsilon}-\varphi_{t}}{\varepsilon}-\frac{d\varphi_{t}}{dt}\big|\big|_{\mathfrak K,\Lam}=
\lim_{\varepsilon\to 0}\big|\big|\sum_{s=1}^\infty(t^s-\frac{(t+\varepsilon)^s-t^s}{\varepsilon s})h^s\big|\big|_{\mathfrak K,\Lam}$

$\le \lim_{\varepsilon\to 0}\sum_{s=1}^\infty\big|t^s-\frac{(t+\varepsilon)^s-t^s}{\varepsilon s}\big| ||h^s||_{\mathfrak K,\Lam}\le \lim_{\varepsilon\to 0}\sum_{s=1}^\infty\big|t^s-\frac{(t+\varepsilon)^s-t^s}{\varepsilon s}\big|=0$ par un argument direct de convergence dominée.
\end{preuve}
\begin{coro}\label{orange}
on a $\frac{d}{dt}||\ln(1-th)||^2_{\mathfrak K,\Lam}=\la\frac{d}{dt}\ln(1-th),\ln(1-th)\ra_{\mathfrak K,\Lam}+\la\ln(1-th),\frac{d}{dt}\ln(1-th)\ra_{\mathfrak K,\Lam}.$
\end{coro}
\begin{preuve}{du corollaire \ref{orange}}
Avec les notations du lemme \ref{ananas}, on a :

$\frac{1}{\varepsilon}\left(\la \varphi_{t+\varepsilon},\varphi_{t+\varepsilon}\ra_{\mathfrak K,\Lam}-\la \varphi_{t},\varphi_{t}\ra_{\mathfrak K,\Lam}\right)-
\la \frac{d\varphi_{t}}{dt},\varphi_{t}\ra_{\mathfrak K,\Lam}-\la \varphi_{t}, \frac{d\varphi_{t}}{dt}\ra_{\mathfrak K,\Lam}=$

$
\la \varphi_{t+\varepsilon},\frac{\varphi_{t+\varepsilon}-\varphi_{t}}{\varepsilon}\ra_{\mathfrak K,\Lam}+\la \frac{\varphi_{t+\varepsilon}-\varphi_{t}}{\varepsilon},\varphi_{t}\ra_{\mathfrak K,\Lam}-
\la \frac{d\varphi_{t}}{dt},\varphi_{t}\ra_{\mathfrak K,\Lam}-\la \varphi_{t}, \frac{d\varphi_{t}}{dt}\ra_{\mathfrak K,\Lam}
=
\la \varphi_{t+\varepsilon}-\varphi_{t},\frac{\varphi_{t+\varepsilon}-\varphi_{t}}{\varepsilon}\ra_{\mathfrak K,\Lam}+$

$\la \varphi_{t},
\frac{\varphi_{t+\varepsilon}-\varphi_{t}}{\varepsilon}-\frac{d\varphi_{t}}{dt}\ra_{\mathfrak K,\Lam}+
\la\frac{\varphi_{t+\varepsilon}-\varphi_{t}}{\varepsilon}-\frac{d\varphi_{t}}{dt},\varphi_{t}\ra_{\mathfrak K,\Lam}$. 

La démonstration s'achève avec l'inégalité de Cauchy-Schwarz associée au lemme \ref{ananas}.
\end{preuve}

 \emph{Quatrième  point}: On démontre l'équation(\ref{bisbis}).

on a $||\ln(1-th)||^2_{\mathfrak K,\Lam}=-2\big(\mathfrak{S}_{1}\mu_{1}(f_{t})+\mathfrak{S}_{2}\mu_{2}(f_{t})\big)$. D'où en utilisant le corollaire \ref{orange} :

$-2\frac{d}{dt}\big(\mathfrak{S}_{1}\mu_{1}(f_{t})+\mathfrak{S}_{2}\mu_{2}(f_{t})\big)=\la \frac{d\ln(1-th)}{dt},\ln(1-th)\ra_{\mathfrak K,\Lam}+\la \ln(1-th),\frac{d\ln(1-th)}{dt}\ra_{\mathfrak K,\Lam}.$
 Or \begin{equation}\label{zip}
\widehat{ \frac{d}{dt}\ln f_{t}}(u,v)=\int_{\tore^2}\frac{d}{dt}\ln f_{t}\;\;\chi_{2}^{-u}\chi_{2}^{-v}d\sigma_{2}=\int_{\tore^2}\frac{-h}{1-th}\chi_{1}^{-u}\chi_{2}^{-v}d\sigma_{2}=\frac{1}{t}\widehat{\frac{1}{f_{t}}}(u,v),
 \end{equation}
 D'où
  
 $\la \frac{d\ln(1-th)}{dt},\ln(1-th)\ra_{\mathfrak K,\Lam}=\frac{1}{t}\sum_{\zz^2}(\mathfrak{S}_{1}|u|+ \mathfrak{S}_{2}|v|) \widehat{\frac{1}{f_{t}}}(u,v)\overline{\widehat{\ln f_{t}}(u,v)}=-\frac{1}{t}\big(\mathfrak{S}_{1}c_{1}(f_{t})+ \mathfrak{S}_{2}c_{2}(f_{t})\big)$, d'après les notations \ref{couic}.
 
 Par ailleurs, 
 
 $\sum_{\zz^2}|u|\;\widehat{\ln f_{t}}(u,v)\overline{\widehat{\frac{d}{dt}\\ln f_{t}}(u,v)}$
 $=\sum_{\zz^2}|u|\;\widehat{\ln f_{t}}(u,v)\widehat{\frac{d}{dt}\ln f_{t}}(-u,-v)$ (car $ln f_{t}$ est réelle)
 
 $\sum_{\zz^2}|u|\;\widehat{\ln f_{t}}(u,v)\frac{1}{t}\widehat{ \frac{1}{f_{t}}}(-u,-v)$ d'après l'égalité (\ref{zip})
 
 $=\sum_{\zz^2}|-u|\;\overline{\widehat{\ln f_{t}}(-u,-v)}\frac{1}{t}\widehat{ f_{t}}(-u,-v)=\frac{1}{t}\sum_{\zz^2}|u|\;\overline{\widehat{\ln f_{t}}(u,v)}{\widehat{\frac{1}{ f_{t}}}(u,v)}$. 
 On en déduit que 
 
 $\frac{d}{dt}\big(\mathfrak{S}_{1}\mu_{1}(f_{t})+\mathfrak{S}_{2}\mu_{2}(f_{t})\big)=\frac{1}{t}\big(\mathfrak{S}_{1}c_{1}(f_{t})+ \mathfrak{S}_{2}c_{2}(f_{t})\big)$ et par intégration l'égalité (\ref{bisbis}).
 
  \bibliography{Toeplitz}
\end{document}